# THE LINEAGE PROCESS IN GALTON–WATSON TREES AND GLOBALLY CENTERED DISCRETE SNAKES

By Jean-François Marckert

*Université Bordeaux*

We consider branching random walks built on Galton–Watson trees with offspring distribution having a bounded support, conditioned to have $n$ nodes, and their rescaled convergences to the Brownian snake. We exhibit a notion of "globally centered discrete snake" that extends the usual settings in which the displacements are supposed centered. We show that under some additional moment conditions, when $n$ goes to $+\infty$, "globally centered discrete snakes" converge to the Brownian snake. The proof relies on a precise study of the lineage of the nodes in a Galton–Watson tree conditioned by the size, and their links with a multinomial process [the lineage of a node $u$ is the vector indexed by $(k,j)$ giving the number of ancestors of $u$ having $k$ children and for which $u$ is a descendant of the $j$th one]. Some consequences concerning Galton–Watson trees conditioned by the size are also derived.

## 1. Introduction.

1.1. *A model of centered discrete snake.* We first begin with the formal description of the notion of trees and branching random walks.

Let $\mathbb{U}=\{\varnothing\}\cup\bigcup_{n\geq 1}\mathbb{N}^{\star n}$ be the set of finite words on the alphabet $\mathbb{N}^\star=\{1,2,\ldots\}$. For $u=u_1\ldots u_n$ and $v=v_1\ldots v_m\in\mathbb{U}$, we let $uv=u_1\ldots u_nv_1\ldots v_m$ be the concatenation of the words $u$ and $v$ (by convention, $\varnothing u=u\varnothing=u$). Following Neveu [22], we call planar tree $T$ a subset of $\mathbb{U}$ containing the root $\varnothing$, and such that if $ui\in T$ for some $u\in\mathbb{U}$ and $i\in\mathbb{N}^\star$, then $u\in T$ and for all $j\in[\![1,i]\!]$, $uj\in T$. The elements of a tree are called nodes or vertices. For $i\neq j$, the nodes $ui$ and $uj$ are called brothers and $u$ their father. We let $c_u(T)=\sup\{i:ui\in T\}$ be the number of children of $u$ [here $c_u(T)$ will be always finite]. A node without any child is called a leaf, and we denote by $\partial T$ the set of leaves of $T$. If $v \neq \varnothing$, we say that $uv$ is a descendant of $u$ and $u$ is an ancestor of $uv$. An edge is a pair $\{u, v\}$ where $u$ is the father of $v$. A path $[\![u, v]\!]$ between the nodes $u$ and $v$ in a tree $T$ is the (minimal) sequence of nodes $u := u_0, \dots, u_j := v$ such that, for any $i \in [\![0, j-1]\!]$, $\{u_i, u_{i+1}\}$ is an edge. Set also $]\!]u, v[\![ = [\![u, v]\!] \setminus \{u, v\}$ and similar notation for $[\![u, v[\![$ and for $]\!]u, v]\!]$. The distance $d_T$, or simply $d$, is the usual graph distance. The depth of $u$ is $|u| = d(\varnothing, u)$. The cardinality of $T$ is denoted by $|T|$, and we let $\mathcal{T}$ (resp. $\mathcal{T}_n$) be the set of planar trees (resp. with $n$ edges, i.e., $n+1$ vertices).

A *branching walk* is a pair $(T, \ell)$ where $T$ is a tree—called the underlying tree—and $\ell$, the label function, is an application from $T$ taking its values in $\mathbb{R}$. In other words, it is a tree in which every vertex owns a real label. We let $\mathcal{B}$ be the set of branching walks, and $\mathcal{B}_n$ be the branching walks associated with trees from $\mathcal{T}_n$.

We introduce now some randomness and construct a probability distribution on $\mathcal{B}$ and on $\mathcal{B}_n$.

The set of underlying trees is endowed with the distribution of the family tree of a Galton–Watson (GW) process with offspring distribution $\boldsymbol{\mu} = (\mu_k)_{k \geq 0}$ starting from one individual. In this model, all the nodes have a random number of children, according to the distribution $\mu$, independently from the other individuals. We denote by $\mathbf{T}$ a random tree under this distribution (see, e.g., [1, 10] and most of the cited papers for more information on GW processes and trees).

The distribution of the labels is defined as follows. Consider $(\nu_k)_{k \in \{1,2,\dots\}}$ a family of distributions, where $\nu_k$ is a distribution on $\mathbb{R}^k$. The labels are defined conditionally on the underlying tree $\mathbf{T}$: Set $\ell(\varnothing) = 0$, and for any $u \in \mathbf{T} \setminus \partial \mathbf{T}$, consider

$$X_u := (\ell(u1) - \ell(u), \dots, \ell(uc_u(\mathbf{T})) - \ell(u)),$$

the evolution-vector of the labels between $u$ and its children. Conditionally on $\mathbf{T}$, we assume that the r.v. $X_u$ are independent, and that $X_u$ has distribution $\nu_{c_u(\mathbf{T})}$. This determines a distribution on $\mathcal{B}$. For example, if $\nu_k$ is the uniform distribution on $\{-1, +1\}^k$ for any $k > 0$, then the r.v. $\ell(u1) - \ell(u), \dots, \ell(uc_u(\mathbf{T})) - \ell(u)$ are independent with common distribution $\frac{1}{2}(\delta_{+1} + \delta_{-1})$ ($\delta_x$ stands for the Dirac mass at $x$). In the case where $\nu_k$ is the uniform distribution on $\{(1, \dots, k), (-1, \dots, -k)\}$, the r.v. $\ell(ui) - \ell(u)$ and $\ell(uj) - \ell(u)$ are not independent and do not have the same distribution.

Notice that a sequence of i.i.d. $\mu$-distributed random variables indexed by $\mathbb{U}$ allows to build the Galton–Watson trees, and a sequence of random variables indexed by $\mathbb{U} \times \mathbb{N}$ allows to build all the labels (by attaching to the elements of $\mathbb{U}$ a list of random variables with distribution $\nu_1, \nu_2, \dots$). We assume that we work on an underlying probability space $(\Omega, \mathcal{A}, \mathbb{P})$ on which are defined all the random variables and processes used in this paper.

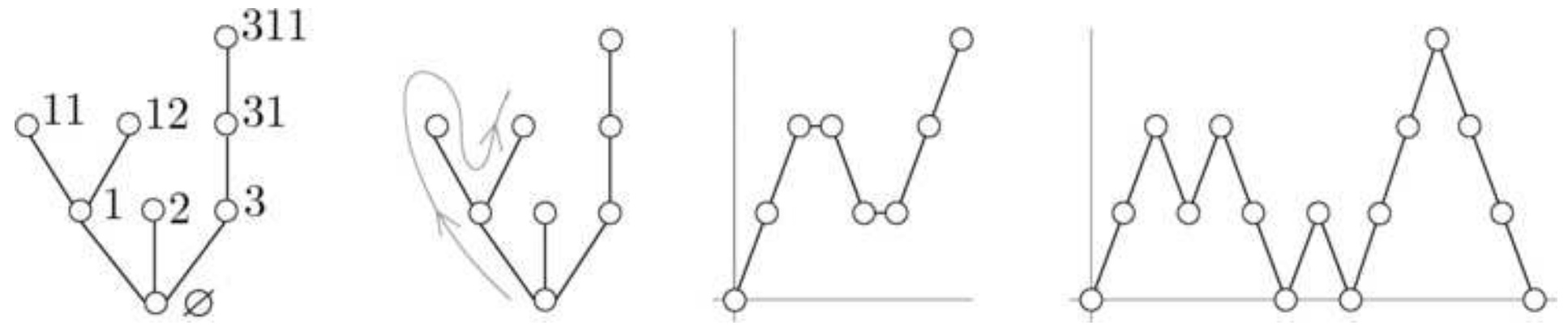

FIG. 1. *A tree on which is indicated the depth-first traversal, its height and contour processes.*

We define now two sets of assumptions $(\mathrm{H}_1)$ and $(\mathrm{H}_2)$ that will be assumed to be satisfied in most of our results. $(\mathrm{H}_1)$ is the following set of conditions: $\mu$ is nondegenerate, critical and has a bounded support,

$$(\mathrm{H}_1) := \left(\mu_0 + \mu_1 \neq 1, \sum_{k\geq 0} k\mu_k = 1, \text{ there exists } K > 0 \text{ s.t. } \sum_{k\leq K} \mu_k = 1\right).$$

Under $(\mathrm{H}_1)$ the variance $\sigma_\mu^2$ of $\mu$ is finite and nonzero. The bounded support condition is quite a strong restriction, but considering nonbounded distribution leads to nontrivial complications, and we were unable to extend to that case the most important results. The condition on the mean can be seen as a normalization, since any distribution $\tilde{\mu}$ related to $\mu$ by $\tilde{\mu}_k = \mu_k a^k/(\sum_j a^k \mu_k)$ for some $a > 0$ induces the same distribution as $\mu$ on GW-trees conditioned by the size.

Let $Y^{(k)} = (Y_{k,1}, \dots, Y_{k,k})$ be $\nu_k$-distributed. We denote by $\nu_{k,j}$, $m_{k,j}$ and $\sigma_{k,j}^2$ the distribution, the mean and the variance of $Y_{k,j}$. We call *global mean* and *global variance* of the branching random walk,

$$\mathbf{m} = \sum_{k\geq 1}\sum_{j=1}^{k} \mu_k m_{k,j} \quad \text{and} \quad \boldsymbol{\beta}^2 = \sum_{k\geq 1}\sum_{j=1}^{k} \mu_k \mathbb{E}(Y_{k,j}^2).$$

Let $(\mathrm{H}_2)$ denote the conditions that the global mean is null, the global variance finite and positive, and for a $p > 4$, the centered $p$th moment of the $Y_{k,j}$'s is finite:

$$(\mathrm{H}_2) := \begin{pmatrix} \mathbf{m} = 0 \text{ and } \beta \in (0, +\infty), \text{ there exist } p > 4 \text{ s.t. for} \\ \text{any } (k,j), 1 \leq j \leq k \leq K, \mathbb{E}(|Y_{k,j} - m_{k,j}|^p) < +\infty. \end{pmatrix}.$$

*Encoding of branching random walks.* We denote by $\preccurlyeq$ the lexicographical order (LO) on the planar trees (and $u \prec v$ if $u \preccurlyeq v$ and $u \neq v$), and let $u(k)$ be the $k$th vertex in the LO $[u(0) = \varnothing]$.

We study the asymptotic behavior of branching random walks via their encoding by depth-first-traversal. The depth-first traversal of a tree $T \in \mathcal{T}_n$ is a function

$$F_T : \{0, \dots, 2n\} \to \{\text{vertices of } T\},$$

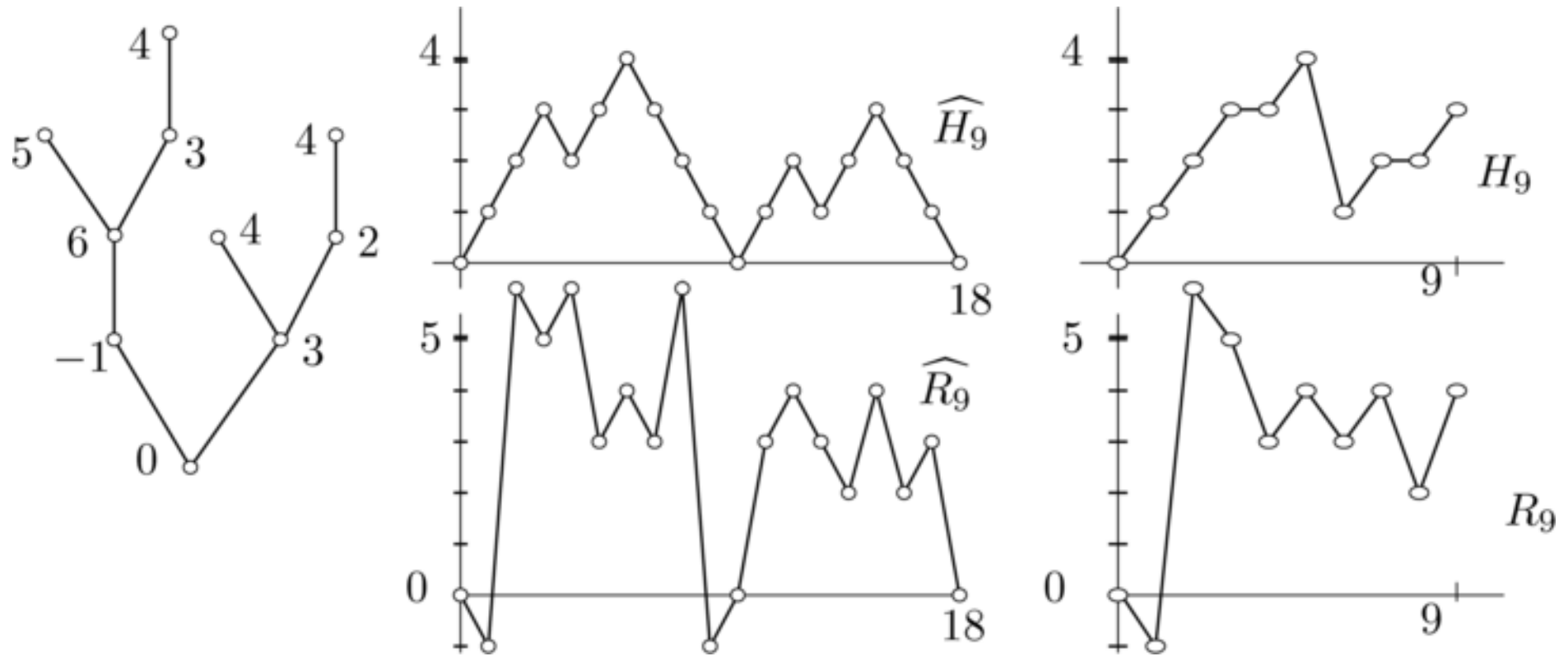

Fig. 2. *A branching random walk from $\mathcal{B}_9$. On the first column, the contour process and the contour label process, on the second column, the height process and the height label process.*

which we regard as a walk around $T$, as follows: $F_T(0) = \varnothing$, and given $F_T(i) = z$, choose if possible and according to the LO, the smallest child $w$ of $z$ which has not already been visited, and set $F_T(i+1) = w$. If not possible, let $F_T(i+1)$ be the father of $z$.

We now encode the branching random walk with the help of a pair of processes. For any $k \in [\![0, |T|-1]\!]$, let $H_k^T = |u(k)|$ and $R_k^T = \ell(u(k))$. The height process $(H_s^T, s \in [0, |T|-1])$ and label process $(R_s^T, s \in [0, |T|-1])$ are obtained from the sequences $(H_k^T)$ and $(R_k^T)$ by linear interpolation. Alternatively, one may encode the branching random walk with a pair of processes associated with the depth-first traversal: for any $k \in [\![0, 2(|T|-1)]\!]$, let $\widehat{H}_k^T = |F_T(k)|$ and $\widehat{R}_k^T = \ell(F_T(k))$. The processes $(\widehat{H}_s^T, s \in [0, 2(|T|-1)])$ and $(\widehat{R}_s^T, s \in [0, 2(|T|-1)])$, obtained by interpolation, are called respectively the contour process and the contour label process; the pair $(\widehat{H}^T, \widehat{R}^T)$ is called the head of the discrete snake. See some illustrations on Figures 1 and 2.

Let $d := \gcd\{k, k \geq 1, \mu_k > 0\}$. The support of the distribution of $|\mathbf{T}|$—we write $\mathrm{supp}(|\mathbf{T}|)$—is included in $1 + d\mathbb{N}$ [and $\mathbb{P}(|\mathbf{T}| = 1 + kd) > 0$ for every $k$ large enough]. For $n + 1 \in \mathrm{supp}(|\mathbf{T}|)$, the distribution $\mathbb{P}$ under the conditioning by $|\mathbf{T}| = n + 1$ is denoted by $\mathbb{P}_n$, in other words

$$\mathbb{P}_n = \mathbb{P}(\cdot \,|\, |\mathbf{T}| = n + 1).$$

Even if not recalled, each statement concerning weak convergence under $\mathbb{P}_n$ is assumed to be along the subsequence $(n_k)_k$ for which $\mathbb{P}_{n_k}$ is well defined. In the proofs we will treat only the case $d = 1$, the general case being treated with slight modifications.

We define $\mathbf{h}_n$, $\widehat{\mathbf{h}}_n$, $\mathbf{r}_n$ and $\widehat{\mathbf{r}}_n$ to be the processes $H^{\mathbf{T}}$, $\widehat{H}^{\mathbf{T}}, R^{\mathbf{T}}$ and $\widehat{R}^{\mathbf{T}}$ under $\mathbb{P}_n$, interpolated as follows:

$$\mathbf{h}_n(s) = \frac{H^{\mathbf{T}}_{ns}}{n^{1/2}}, \qquad \widehat{\mathbf{h}}_n(s) = \frac{\widehat{H}^{\mathbf{T}}_{2ns}}{n^{1/2}},$$

$$\mathbf{r}_n(s) = \frac{R^{\mathbf{T}}_{ns}}{n^{1/4}}, \qquad \widehat{\mathbf{r}}_n(s) = \frac{\widehat{R}^{\mathbf{T}}_{2ns}}{n^{1/4}} \qquad \text{for any } s \in [0,1].$$

THEOREM 1. *If* $(\mathrm{H}_1)$ *and* $(\mathrm{H}_2)$ *are satisfied, then*

$$(\mathbf{h}_n, \widehat{\mathbf{h}}_n, \mathbf{r}_n, \widehat{\mathbf{r}}_n) \xrightarrow[n]{(d)} (\mathbf{h}, \mathbf{h}, \beta\mathbf{r}, \beta\mathbf{r})$$

*in* $C([0,1],\mathbb{R}^4)$ *endowed with the topology of uniform convergence, where* $\mathbf{h} = 2\mathsf{e}/\sigma_\mu$ *and* $\mathsf{e}$ *is the normalized Brownian excursion, and where, conditionally on* $\mathbf{h}$*,* $\mathbf{r}$ *is a centered Gaussian process with covariance function*

$$\operatorname{cov}(\mathbf{r}(s), \mathbf{r}(t)) = \check{\mathbf{h}}(s,t) := \min_{u \in [s\wedge t, s\vee t]} \mathbf{h}(u) \qquad \textit{for any } s,t \in [0,1].$$

Notice that the same processes $\mathbf{h}$ and $\mathbf{r}$ appear twice in the limit process. The convergence of processes associated with the contour processes (with a $\widehat{\ }$ ) to the same limit as the one associated with the height processes is well understood now, and "almost" generic (Duquesne and Le Gall [9], Section 2.5, and [21]), we then concentrate only on the height process. The process $(\mathbf{r}, \mathbf{h})$ (or with a different scaling) will be called the *head of the Brownian snake with lifetime process the normalized Brownian excursion* (*BSBE*). We refer to the works of Le Gall (e.g., [16] and with Duquesne [10]) for information on the Brownian snake and to the papers cited below for discrete approaches to this object.

In the present work we deal only with the head of the snakes; this is, in principle, different than snakes even if, thanks to the homeomorphism theorem [20] evoked below, Theorem 1 has some direct interpretation in terms of snakes. We refer to [13, 20] for the notion of discrete snake which is the discrete analogue of BSBE: the discrete snake associated with the branching random walk $(T, \ell)$ is the pair $(\widehat{H}^T, \Phi)$ where $\Phi = (\Phi_k)_{k \in [\![0, 2(|T|-1)]\!]}$ and $\Phi_k$ is the sequence of labels on the branch $[\![\varnothing, F_T(k)]\!]$.

*Related works.* The convergence $\widehat{\mathbf{h}}_n \xrightarrow[n]{(d)} \mathbf{h}$ is due to Aldous [1, 2] (see also Marckert and Mokkadem [21] for a revisited proof, Pitman [25], Chapters 5 and 6 and Duquesne [9] and Duquesne and Le Gall [10], Section 2.5, for generalization to GW trees with offspring distribution having infinite variance).

The two first results concerning the convergence of discrete snakes to the BSBE appeared in two independent works:

- Chassaing and Schaeffer [7] deal with discrete snakes built on underlying trees chosen uniformly in $\mathcal{T}_n$ [this corresponds to the case $\mu \sim \mathit{Geom}(1/2)$] and where the displacements are i.i.d., and for any $k, j$, $\nu_{k,j}$ is the uniform distribution in $\{-1, 0, +1\}$. They show the convergence of the head of the discrete snake for the Skohorod topology, and the convergence of the moments of the maximum of $\mathbf{r}_n$ are also given. This study was motivated by the deep relation between this model of discrete snake and random rooted quadrangulations, underlined by the authors.
- Marckert and Mokkadem [20] studied also the case $\mu \sim \mathit{Geom}(1/2)$, but with more general centered displacements that have moments of order $6+\varepsilon$ (the distribution $\nu_{k,j}$ does not depend on $k, j$, but $\nu_k$ is not assumed to be $\nu_{k,1} \times \cdots \times \nu_{k,k}$). The convergence of the head of the snake holds in $(C[0,1], \mathbb{R}^2)$ and the convergence of the snake itself is given thanks to a "homeomorphism theorem" which implies that the convergence of the snake and of its tour (in space of continuous functions) are equivalent. Here it implies that, under the hypothesis of Theorem 1, the discrete snake associated with our model of labeled trees converges weakly to the BSBE.

Then some generalizations appeared few months later:

- Gittenberger [11] provides a generalization of a lemma from [20] allowing him to consider snakes with underlying GW trees conditioned by the size (condition equivalent to $\mathrm{H}_1$). The displacements must be centered and have moments of order $8+\varepsilon$.
- Janson and Marckert [13] show that, in the i.i.d. case [$\nu_{k,j}$ do not depend on $(k, j)$], moments of order $4+\varepsilon$ are necessary and sufficient to get the convergence to the BSBE. If no such moment exists, the convergence to a "hairy snake" is proved under the Hausdorff topology.
- In Marckert and Miermont [19], the case of $\nu_{k,j}$ depending on $k, j$ is investigated (also the underlying GW trees are allowed to have two types). The hypothesis are for each $k, j$, $m_{k,j} = 0$, condition $(\mathrm{H}_2)$ is satisfied, and then $\sum_{k,j} \mu_k \sigma_{k,j}^2 < +\infty$. A motivation was to generalize the works of Chassaing and Schaeffer [7] concerning quadrangulations to bipartite maps.

Another important point is the convergence of the occupation measure of the head of the discrete snake to the one of the BSBE, the random measure named ISE (the integrated superBrownian excursion introduced by Aldous [3]; see also Le Gall [16] and [13, 20]). Using the convergence of the discrete snake to the BSBE, Bousquet–Mélou [4] and Bousquet–Mélou and Janson [5] deduce new results on the ISE and on the BSBE; for example, some properties on the support of the ISE, and of random density of the ISE are derived. We refer also to Le Gall [15] for the convergence of the discrete snake conditioned to stay positive.

The novelty in the present paper is that the condition $\{m_{k,j}=0, \forall k,j\}$ is replaced by $\mathbf{m}=\sum_{k\geq 1}\sum_{j=1}^{k}\mu_k m_{k,j}=0$. This allows to consider some natural models where, for example, the displacements are not random knowing the underlying tree (see Section 1.3). The proof of Theorem 1 relies in part on some results from [19], and on a new approach necessary to control the contribution of the mean of the displacements; the main point for this is the comparison of the lineage of each node, with some multinomial r.v. This is the aim of Theorem 2, that we think interesting in itself, since it reveals a thin global behavior of GW trees conditioned by the size. Unfortunately, the price of this generalization is to consider only offspring distribution with bounded support. The reason comes from the proof of Theorem 2. We guess that some generalization for all families of GW trees (with finite variance) may be found, but for this a control of an infinite sequence of processes arising in Theorem 2 should be provided for what we were unable to do.

1.2. *On the lineage of nodes.* Assume that $(\mathrm{H}_1)$ and $(\mathrm{H}_2)$ hold, and let $K$ be a bound on the support of the offspring distribution. We work again conditionally on $\mathbf{T}$. For any node $u=i_1\dots i_h\in\mathbf{T}$, let $u_j=i_1\dots i_j$ and $[\![\varnothing,u]\!]=\{\varnothing=u_0,u_1,\dots,u_h\}$ be the ancestral line of $u$ back to the root. Conditionally on $\mathbf{T}$, $\ell(u)$ owns the following representations:

$$\ell(u)=\sum_{m=1}^{|u|}\ell(u_m)-\ell(u_{m-1}), \tag{1}$$

where $\ell(u_m)-\ell(u_{m-1})$ is $\nu_{k,j}$-distributed when $c_{u_{m-1}}(\mathbf{T})=k$ and $i_m=j$, and where the r.v. $(\ell(u_m)-\ell(u_{m-1}))$'s are independent (conditionally on $\mathbf{T}$); the variables $\ell(u_m)-\ell(u_{m-1})$ will be often called displacements.

Consider the array

$$I_K=\{(k,j),1\leq j\leq k\leq K\}.$$

Let $T\in\mathcal{T}$ and $u$ be a node of $T$. For any $(k,j)\in I_K$, let $A_{u,k,j}(T)$ be the number of strict ancestors $v$ of $u$ (the nodes $v\in[\![\varnothing,u[\![$) such that $c_v(T)=k$, and such that $u$ is a descendant of $vj$, the $j$th child of $v$ [we write $f_v(u)=j$]. We say that $v$ is an ancestor of type $k,j$ of $u$, and we call the vector $A_u=(A_{u,i})_{i\in I_K}$ the lineage of $u$ (or the content of $[\![\varnothing,u[\![$). See Figure 3.

By (1), conditionally on $\mathbf{T}$, the label $\ell(u)$ owns the following representations:

$$\ell(u)\stackrel{(d)}{=}\sum_{(k,j)\in I_K}\sum_{l=1}^{A_{u,k,j}(\mathbf{T})}Y_{k,j}^{(l)},$$

where the r.v. $Y_{k,j}^{(l)}$ are independent, and where for any $l$, $Y_{k,j}^{(l)}$ is $\nu_{k,j}$ distributed. In order to make more apparent the contribution of the means

$m_{k,j}$'s, and using that $\mathbf{m}=0$, write

$$\ell(u) \overset{(d)}{=} \sum_{(k,j)\in I_K} \sum_{l=1}^{A_{u,k,j}(\mathbf{T})} (Y_{k,j}^{(l)} - m_{k,j}) + \sum_{(k,j)\in I_K} (A_{u,k,j}(\mathbf{T}) - \mu_k |u|) m_{k,j}. \tag{2}$$

Assume that $\mathbf{T}$ is $\mathbb{P}_n$ distributed, and that $u = u(ns)$ for some $s \in (0,1)$. Conditionally on $|u|$, we will see that both parts of the right-hand side of (2) divided by $n^{1/4}$ converge in distribution, and the limit r.v. are asymptotically independent: in the first part, the fluctuations of $A_{u,k,j}$ around $\mu_k|u|$ are not important when they are crucial in the second sum.

We now concentrate on the r.v. $(A_u)'s$ under $\mathbb{P}_n$. For any $l \in [\![0,n]\!]$, $(k,j) \in I_K$, set

$$\mathbf{g}_{(k,j)}^{(n)}(l) := A_{u(l),k,j} - \mu_k |u(l)|.$$

For every $(k,j) \in I_K$, the process $l \mapsto \mathbf{g}_{(k,j)}^{(n)}(l)$ encodes the evolution of the number of ancestors of type $k,j$ of $u(l)$, when $l$ varies. Consider $\mathbf{G}^{(n)} = (\mathbf{G}^{(n)}(s))_{s\in[0,1]}$ the process taking its values in $\mathbb{R}^{I_K}$ defined by, for any $s$, $\mathbf{G}^{(n)}(s) = (\mathbf{G}_{k,j}^{(n)}(s))_{(k,j)\in I_K}$, where $s \mapsto \mathbf{G}_{k,j}^{(n)}(s)$ is the real continuous process that interpolates $\mathbf{g}_{k,j}^{(n)}$ as follows:

$$\mathbf{G}_{k,j}^{(n)}(s) := \frac{\mathbf{g}_{k,j}^{(n)}(\lfloor ns \rfloor) + \{ns\}(\mathbf{g}_{k,j}^{(n)}(\lfloor ns+1 \rfloor) - \mathbf{g}_{k,j}^{(n)}(\lfloor ns \rfloor))}{n^{1/4}}, \qquad s \in [0,1], \tag{3}$$

where $\{x\}$ stands for the rational part of $x$. The random process $\mathbf{G}^{(n)}$ encodes the lineage of all the nodes of $\mathbf{T}$; its limiting behavior is described by the following theorem.

THEOREM 2. *Under* $(\mathrm{H}_1)$ *and* $(\mathrm{H}_2)$ *the following convergence in distribution holds in* $C([0,1])^{\#I_K} \times C[0,1]$ *endowed with the topology of the uniform convergence*

$$(\mathbf{G}^{(n)}, \mathbf{h}_n) \xrightarrow[n]{(d)} (\mathbf{G}, \mathbf{h}),$$

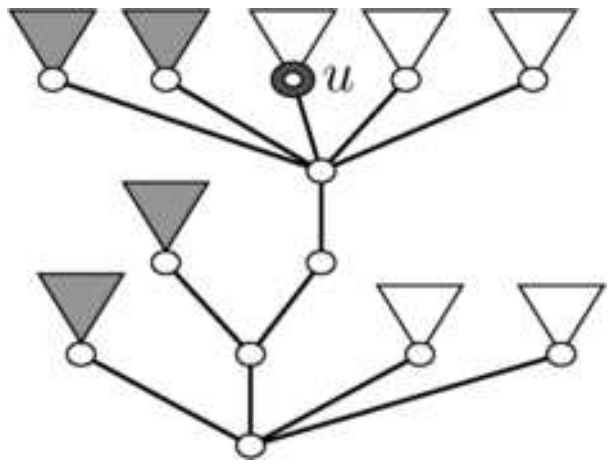

FIG. 3. *On this tree* $A_{u,1,1}=1, A_{u,2,2}=1, A_{u,4,2}=1, A_{u,5,3}=1$, *the others* $A_{u,i}$ *are 0.*

*where* $\mathbf{h}$ *is defined as in Theorem 1, and where conditionally on* $\mathbf{h}$, $\mathbf{G} = (\mathbf{G}_{k,j}(s))_{(k,j)\in I_K, s\in[0,1]}$ *is a real centered Gaussian field with the following covariance function: for any* $(k,j)$ *and* $(k',j')$ *in* $I_K$, $s$ *and* $s'$ *in* $[0,1]$,

$$\operatorname{cov}(\mathbf{G}_{k,j}(s), \mathbf{G}_{k',j'}(s')) = (-\mu_k\mu_{k'} + \mu_k \mathbb{1}_{(k,j)=(k',j')})\check{\mathbf{h}}(s,s'). \tag{4}$$

1.3. *Comments, examples and applications.* (1) Theorem 2 may be considered as the strongest result of this paper. It gives very precise information on the asymptotic behavior of the process $\mathbf{G}^{(n)}$ that encodes the lineage of all the nodes. This gives a "global asymptotic" property reminiscent of the properties of the distinguished branch in "a size biased GW tree" (see [17], Chapter 11).

(2) For any fixed $(k,j) \in I_K$, knowing $\mathbf{h}$, $\mathbf{G}_{k,j}$ is a Gaussian process with covariance function

$$\operatorname{cov}(\mathbf{G}_{k,j}(s), \mathbf{G}_{k,j}(s')) = (-\mu_k^2 + \mu_k)\check{\mathbf{h}}(s,s').$$

In other words, the process $(\mathbf{G}_{k,j}, \mathbf{h})$ has the same distribution as $(\sqrt{-\mu_k^2+\mu_k}\mathbf{r}, \mathbf{h})$, and then up to some multiplicative constants, $(\mathbf{G}_{k,j}, \mathbf{h})$ is the head of a BSBE. As a simple consequence of Theorem 2, we have that $(\mathbf{G}_{k,j}, \mathbf{h})_{(k,j)\in I_K}$ is a sequence of heads of BSBE, and that for any $(k,j) \in I_K$,

$$(\mathbf{G}_{k,j}^{(n)}, \mathbf{h}_n) \xrightarrow[n]{(d)} (\mathbf{G}_{k,j}, \mathbf{h}). \tag{5}$$

The dependence between the different processes $\mathbf{G}_{k,j}$ is ruled out by (4). For any families of real numbers $(\lambda_{k,j})_{(k,j)\in I_K}$, we have

$$\left(\sum_{(k,j)\in I_K} \lambda_{k,j}\mathbf{G}_{k,j}^{(n)}, \mathbf{h}_n\right) \xrightarrow[n]{(d)} \left(\sum_{k,j} \lambda_{k,j}\mathbf{G}_{k,j}, \mathbf{h}\right). \tag{6}$$

(3) Consider the case $\mu = \frac{1}{2}(\delta_0 + \delta_2)$, $\nu_2 = \delta_{(+1,-1)}$, of binary trees in which the displacements are not random: $\ell(u1) - \ell(u) = +1$ and $\ell(u2) - \ell(u) = -1$. We have $\mathbf{m} = 0$ and $\beta^2 = \frac{1}{2}(1+1) = 1$ and Theorems 1 and 2 apply. Hence, the clear positive bias for $R_n(t)$ for small values of $t$ disappears at the limit. Note also that this normalizing factor is exactly the same as if $\nu_2 = \frac{1}{2}(\delta_{(+1,-1)} + \delta_{(-1,+1)})$ [case where $(\ell(u1) - \ell(u), \ell(u2) - \ell(u))$ is equally likely $(+1,-1)$ or $(-1,+1)$] and as if $\nu_2 = (\frac{1}{2}(\delta_{+1} + \delta_{-1}))^2$ [case where the $\ell(u1) - \ell(u)$ and $\ell(u2) - \ell(u)$ are i.i.d., uniform on $\{-1,1\}$]. The question of the convergence of the discrete snake in the case $\nu_2 = \delta_{(+1,-1)}$ appears first in Marckert [18] in relation with some properties of the rotation correspondence, and the difference between left and right depth in binary trees. The convergence of $(\mathbf{r}_n)$ is not given in [18], but the convergence of the occupation measure of $\mathbf{r}_n$, "the discrete ISE," to ISE is established. We refer also to Janson [12] for recent developments concerning the same question.

Further, notice that in this model, the label $\ell(u)$ of a vertex $u$ is $\ell(u) = A_{u,2,1} - A_{u,2,2}$, that is, the number of left steps minus the number of right steps necessary to climb from the root to $u$ in the binary tree. The convergence of $(\mathbf{r}_n)$ can be seen directly via the one of $(\mathbf{G}^{(n)})$:

$$(\mathbf{G}_{2,1}^{(n)}, \mathbf{G}_{2,2}^{(n)}, \mathbf{h}_n) \xrightarrow[n]{(d)} (\mathbf{G}_{2,1}, \mathbf{G}_{2,2}, \mathbf{h}), \tag{7}$$

and then $\mathbf{r}_n = \mathbf{G}_{2,1}^{(n)} - \mathbf{G}_{2,2}^{(n)} \xrightarrow[n]{(d)} \mathbf{G}_{2,1} - \mathbf{G}_{2,2}$ which is, conditionally to $\mathbf{h}$ and according to (4), a centered Gaussian process with covariance function $\check{\mathbf{h}}(s,t)$. Here, the convergence of $(\mathbf{r}_n)$ appears to be a consequence of the convergence of $\mathbf{G}_{2,1}$ and $\mathbf{G}_{2,2}$, encoding the right depth and the left depth in binary trees.

We would like to stress on the following point: discrete snake are usually constructed with "two levels of randomness": the underlying trees are random and so are the displacements given the underlying tree, and then BSBE appears to be a natural limit of these objects. Here, we provide some objects with only "one level of randomness" that converges to the Brownian snake. The BSBE appears as a kind of internal complexity measure in trees measuring the difference between the number of ancestors of type $k, j$ and some expected quantities.

## 2. Proofs.

The proofs rely on a precise study of the lineage of the nodes under $\mathbb{P}_n$ and, in particular, on the comparison of $A_u$ with a multinomial random variable. For this reason, we first give some elements on multinomial distributions and on their asymptotic behaviors. We then proceed to the proof of Theorem 2, showing first the convergence of the uni-dimensional distribution, then the convergence of the finite-dimensional distribution. The proof of Theorem 1 is given afterward. We think that some points of view, especially in the description of the distribution of the lineages in trees under $\mathbb{P}_n$, should provide some new approaches to study the trees under $\mathbb{P}_n$.

### 2.1. *Prerequisite on multinomial distributions.*

The contents of this section are quite classical. Consider $\mathbf{p} = (p_i)_{i \in I_K}$ the distribution on $I_K$, defined by

$$p_{k,j} := \mu_k \qquad \text{for any } (k,j) \in I_K.$$

For any $h \geq 1$, let $\mathbb{N}^I[h]$ be the set of elements $c = (c_i)_{i \in I_K}$ of $\mathbb{N}^{\# I_K}$, such that $\sum_{i \in I_K} c_i = h$. We say that $\mathcal{M}^{(h)}$ is a multinomial r.v. with parameter $h$ and $\mathbf{p}$, if, for any $\mathsf{m} = (\mathsf{m}_i)_{i \in I_K}$,

$$\mathbb{Q}_h(\{\mathsf{m}\}) := \mathbb{P}(\mathcal{M}^{(h)} = \mathsf{m}) = \binom{h}{(\mathsf{m}_i)_{i \in I_K}} \prod_{i \in I_K} p_i^{\mathsf{m}_i} \mathbb{1}_{\mathbb{N}^I[h]}(\mathsf{m}),$$

where $\binom{h}{(\mathsf{m}_i)_{i\in I_K}} = h!/(\prod_{i\in I_K} \mathsf{m}_i!)$. Recall that for any $i \in I_K$, $\mathcal{M}_i^{(h)}$ is a binomial r.v. with parameters $n$ and $p_i$.

In order to fit with further considerations, we introduce the $\#I_K$ dimensional real vector $\mathcal{G}(n,h) = (\mathcal{G}_i(n,h))_{i\in I_K}$ defined by

$$\mathcal{G}_{k,j}(n,h) = n^{-1/4}(\mathcal{M}_{k,j}^{(h)} - \mu_k h) \qquad \text{for any } (k,j) \in I_K.$$

Let $\mathcal{G}_\infty = (\mathcal{G}_{\infty,i})_{i\in I_K}$ be a centered Gaussian vector having as covariance function

$$\text{cov}(\mathcal{G}_{\infty,i}, \mathcal{G}_{\infty,i'}) = -p_i p_{i'} + p_i \mathbb{1}_{i=i'} \qquad \text{for any } i,i' \in I_K. \tag{8}$$

PROPOSITION 3. *Let $(h(n))$ be a sequence of positive integers s.t. $h(n)/\sqrt{n} \to \lambda \in (0,+\infty)$. Under* (H$_1$), *we have* $\mathcal{G}(n,h(n)) \xrightarrow[n]{(d)} \sqrt{\lambda}\mathcal{G}_\infty$ *in* $\mathbb{R}^{\#I_K}$.

PROOF. This may be proved using classical tools. As pointed out by E. Rio in a personal discussion, this is also a consequence of the convergence of the empirical process to the Brownian bridge. We only sketch the proof (for $\lambda = 1$): let $(U_l)_l$ be a sequence of i.i.d. r.v. uniform on $[0,1]$. Let $F_n$ be the associated empirical distribution function and $F$ the distribution function of $U$. Denote by $g_n = F_n - F$.

According to Donsker [8], $\sqrt{n}g_n \xrightarrow[n]{(d)} \mathsf{b}$, where $\mathsf{b}$ is a normalized Brownian bridge.

Take $\mathbf{q} = (q_l)_{l\in\mathbb{N}}$ a distribution on $\mathbb{N}$ and consider

$$\mathcal{N}_k^{(n)} = \#\{j, j \in \{1,\dots,n\}, U_j \in [q_1 + \cdots + q_k, q_1 + \cdots + q_{k+1}]\}.$$

Then $(\mathcal{N}_k^{(n)})_{k\ge 1}$ is a multinomial r.v. with parameters $n$ and $\mathbf{q}$ and satisfies

$$(\mathcal{N}_k^{(n)} - q_k n)/\sqrt{n} = \sqrt{n}(g_n(q_1 + \cdots + q_{k+1}) - g_n(q_1 + \cdots + q_k)).$$

By Donsker, for any $L > 0$, $((\mathcal{N}_k^{(n)} - q_k n)/\sqrt{n})_{k\le L}$ converges in distribution to $(\mathsf{b}_{q_1+\cdots+q_{k+1}} - \mathsf{b}_{q_1+\cdots+q_k})_{k\le L}$. The properties of $\mathsf{b}$ allow to conclude. □

The following proposition will be used in the proof of the tightness of $(\mathbf{G}^{(n)})$.

PROPOSITION 4. *Under* (H$_1$), *for any $\beta > 1$, there exists $c > 0$ such that, for any $h > 0$, any $n > 0$,*

$$\mathbb{E}(\|\mathcal{G}(n,h)\|_1^\beta) \le c(h/\sqrt{n})^{\beta/2}.$$

Recall that all the norms are equivalent in $\mathbb{R}^{\# I_K}$. Here, we use $\|X\|_1 = \sum_{(k,j)\in I_K} |X_{k,j}|$.

Proof of Proposition 4. First, since $\|X\|_1^\beta \le c\sum |X_{k,j}|^\beta$ for some $c > 0$,

$$\mathbb{E}(\|\mathcal{G}(n,h)\|_1^\beta) \le c \sum_{(k,j)\in I_K} \mathbb{E}(|n^{-1/4}(\mathcal{M}_{k,j}^{(h)} - \mu_k h)|^\beta).$$

Since $\mathcal{M}_{k,j}^{(h)}$ is a binomial random variable with parameter $\mu_k$ and $h$, $\mathbb{E}(|(\mathcal{M}_{k,j}^{(h)} - \mu_k h)|^\beta) \le C(\mu_k,\beta) h^{\beta/2}$, where the constant $C(\mu_k,\beta)$ depends on $\mu_k$ and $\beta$ (see Petrov [24], Theorem 2.10, page 62). □

2.2. *Decomposition of trees using the lineages.* A forest is a finite sequence of trees, that is an element of $\mathcal{F} := \bigcup_{k\ge 0} \mathcal{T}^k$. For any $k \in \mathbb{N}$, a forest with $k$ roots is a $k$-tuple of planar trees $f = (t^1,\dots,t^k)$. The size $|f|$ of $f$ is $|t^1| + \dots + |t^k|$. We denote by $\mathbf{f}_k = (\mathbf{T}^1,\dots,\mathbf{T}^k)$ a random forest in which the trees $\mathbf{T}^1,\dots,\mathbf{T}^k$ are i.i.d. GW trees with offspring distribution $\mu$. For any $\mathsf{a} = (\mathsf{a}_{k,j})_{(k,j)\in I_K} \in \mathbb{R}^{I_k}$, write

$$N_1(\mathsf{a}) = \sum_{(k,j)\in I_K} (j-1)\mathsf{a}_{k,j} \quad \text{and} \quad N_2(\mathsf{a}) = \sum_{(k,j)\in I_K} (k-j)\mathsf{a}_{k,j}.$$

Proposition 5. *Let $h$ be a nonnegative integer. For any $\mathsf{a} \in \mathbb{N}^I[h]$, and any $m \in [\![0,n]\!]$,*

$$\mathbb{P}_n(A_{u(m)} = \mathsf{a}) = \mathbb{Q}_h(\mathsf{a}) \frac{\mathbb{P}(|\mathbf{f}_{N_1(\mathsf{a})}| = m-h, |\mathbf{f}'_{1+N_2(\mathsf{a})}| = n+1-m)}{\mathbb{P}(|\mathbf{T}| = n)}, \tag{9}$$

*where $\mathbf{f}$ and $\mathbf{f}'$ are two independent forests.*

Proof. To build a tree $T$ of $\mathcal{T}_n$ such that $A_{u(m)} = \mathsf{a}$, we first build the branch $b = [\![\varnothing, u(m)]\!]$: Exactly $\mathsf{a}_{k,j}$ ancestors $v$ among the $h$ strict ancestors of $u$ satisfy $(c_v(T), f_v(u)) = (k,j)$. Hence, there are $\binom{n}{\mathsf{a}}$ ways to build $b$. Then, we complete $b$ in grafting on its neighbors some subtrees satisfying the following constraints. When $A_{u(m)} = \mathsf{a}$, the number of subtrees rooted on the neighbors of the branch $[\![\varnothing, u(m)[\![$ visited before $u(m)$ [resp. after $u(m)$] are respectively

$$N_1(\mathsf{a}) = \#\{w, d([\![\varnothing, u(m)[\![, w) = 1, w \prec u(m)\},$$
$$1 + N_2(\mathsf{a}) = \#(\{u(m)\} \cup \{w, d(]\!]u(m), \varnothing]\!], w) = 1, u(m) \prec w\}).$$

See an illustration on Figure 4. The $N_1(\mathsf{a})$ subtrees must contain exactly $m - |u(m)|$ nodes (the nodes, among the $m+1$ first, not on $[\![\varnothing, u(m)]\!]$), and the $1+N_2(\mathsf{a})$ subtrees must contain exactly $n+1-m$ nodes [the nodes visited after $u(m)$, $u(m)$ included]. In other words, we need two forests containing respectively $m-h$ and $n+1-m$ nodes. Hence, using simple considerations on the probability distribution of GW trees, we get the announced result. □

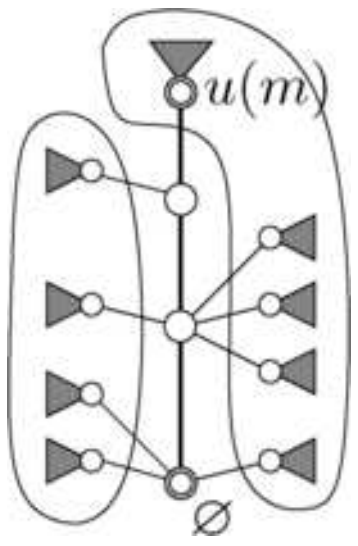

FIG. 4. *The two forests considered in the decomposition.*

A consequence of this proposition is

$$\mathbb{P}_n(|u(m)|=h)$$

$$= \sum_{\mathsf{x}\in\mathbb{N}^I[h]} \mathbb{Q}_h(\mathsf{x}) \frac{\mathbb{P}(|\mathbf{f}_{N_1(\mathsf{x})}|=m-h, |\mathbf{f}'_{1+N_2(\mathsf{x})}|=n+1-m)}{\mathbb{P}(|\mathbf{T}|=n)} \tag{10}$$

$$= \frac{\mathbb{P}(|\mathbf{f}_{N_1(\mathcal{M}^{(h)})}|=m-h, |\mathbf{f}'_{1+N_2(\mathcal{M}^{(h)})}|=n+1-m)}{\mathbb{P}(|\mathbf{T}|=n)}, \tag{11}$$

where $\mathcal{M}^{(h)}$ is a multinomial random variable with parameters $h$ and $\mathbf{p}$.

2.2.1. *Few facts concerning random forests and random trees.* Let $(W_i)_{i\geq 0}$ be a random walk starting from 0 with i.i.d. increments with distribution $(\tilde{\mu}_k)_{k\geq -1} = (\mu_{k+1})_{k\geq -1}$ (i.e., with increment $\xi-1$, where $\xi$ is $\mu$-distributed). We have the following:

LEMMA 6. *Assume* $(\mathrm{H}_1)$ *holds true.*

(i) *(Otter [23]) For any* $k\geq 1$ *and* $n\geq k$, $\mathbb{P}(|\mathbf{f}_k|=n) = \frac{k}{n}\mathbb{P}(W_n=-k)$.
(ii) *[Central local limit theorem (CLLT)]*

$$\sup_{l\in -n+d\mathbb{N}} \left| \frac{\sqrt{n}}{d}\mathbb{P}(W_n=l) - \frac{1}{\sqrt{2\pi}\sigma_\mu}\exp\left(-\frac{l^2}{2\sigma_\mu^2 n}\right)\right| \underset{n}{\to} 0. \tag{12}$$

(iii) $\sup_{n\geq 0}\sup_{x\geq 0} x\mathbb{P}(W_n=x) < +\infty$.

(i) is often called “conjugation of tree principle” or “cyclical lemma,” and may be found in Pitman [25], Chapter 5.1 and is usually attributed to Otter, Kemperman or Dvoretzky-Motzkin.

(ii) is usually called the central local limit theorem (see Breuillard [6] for a state of the art). Recall that $d$ is the span of $\mu$. The support of $W_n$ is included in $-n+d\mathbb{N}=\{u\in\mathbb{Z}, u=-n+di, i\in\mathbb{N}\}$. A consequence of (i) and (ii) is that

$$\mathbb{P}(|\mathbf{T}|=n)\sim\frac{dn^{-3/2}}{\sqrt{2\pi}\sigma_\mu}, \tag{13}$$

the equivalent being taken along the subsequence where the left-hand side is nonnull.

Proof of Lemma 6. (iii) $\sup_{n\ge 0}\sup_{x\ge c\sqrt{n}}\{x\mathbb{P}(W_n=x)\}$ is bounded by the Chebyshev inequality. By (ii), $\sup_{x\le c\sqrt{n}}\sqrt{n}\mathbb{P}(W_n=x)\xrightarrow[n]{}\frac{d}{\sqrt{2\pi}\sigma_\mu}$, then $\sup_{n\ge 0}\sup_{x\le c\sqrt{n}}\sqrt{n}\mathbb{P}(W_n=x)$ is finite. □

The following lemma controls the maximum increment in the process $H$ under $\mathbb{P}_n$.

Lemma 7. *Assume* (H$_1$). *For any $c>0$, there exists $\rho>0$ such that*

$$\mathbb{P}_n\Big(\max_l\{||u(l+1)|-|u(l)||\}\ge\rho\log n\Big)=O(n^{-c}).$$

Proof. The proof deeply relies on the conjugation of the tree principle. Take $n+1$ i.i.d. r.v. $X_1,\dots,X_{n+1}$, $\mu$-distributed. Conditionally on $\sum_{i=1}^{n+1}(X_i-1)=-1$, among the $n+1$ shifted sequences $(X_1,\dots,X_{n+1})$, $(X_2,\dots,X_{n+1},X_1),\dots,(X_{n+1},X_1,\dots,X_n)$, exactly one $(X_1^\star,\dots,X_{n+1}^\star)$ corresponds to a sequence $(c_u,u\in T)$ for a tree $T\in\mathcal{T}_n$ (where the $c_u$ are sorted according to the depth first order), and $(X_1^\star,\dots,X_{n+1}^\star)\stackrel{(d)}{=}(c_u,u\in\mathbf{T})$ for $\mathbf{T}$ under $\mathbb{P}_n$.

The inequality $||u(l)|-|u(l+1)||=h>1$ implies that $|u(l+1)|<|u(l)|$, and the deepest common ancestor $v$ of $u(l+1)$ and $u(l)$ has depth $|u(l+1)|-1$. Assume that the tree is visited according to the reversed LO (the order on the alphabet $\mathbb{N}$ is reversed, but if $z$ is a prefix of $z'$, $z$ is still smaller than $z'$: this amounts to walk around the tree counterclockwise and rank the nodes according to their first visit time). In the reversed LO, the nodes in $]\!]v,u(l)[\![$ are visited consecutively, and each of them has at least one child. Under $\mathbb{P}$, when traversing the tree in the LO (or by symmetry in the reversed LO), the gap between two nodes having zero child is a geometrical r.v. $Geom(\mu_0)$. We work now on the LO order. Denote by $X_1,\dots,X_{n+1}$ i.i.d. random variables $\mu$-distributed and by $G_1,G_2,\dots$ the successive gaps between the zeros:

$$\mathbb{P}\biggl(\max_i G_i \ge \frac{\rho\log n}{2}\biggm|\sum_{i=1}^{n+1}(X_i-1)=-1\biggr) = O\biggl(n^{1/2}\mathbb{P}\biggl(\max_{i\le n} G_i \ge \frac{\rho\log n}{2}\biggr)\biggr)$$
$$= o(n^{-c-1}),$$

for $\rho$ large enough. Note that the first maximum is taken on a random number of terms, a.s. bounded by $n$. By the conjugation of the tree principle, we get the result. $\square$

REMARK 1. Using the same argument, one may control the depth of the last node $u(n)$: for any $c>0$, there exists $\rho>0$ such that

$$\mathbb{P}_n(|u(n)|\ge\rho\log n) = O(n^{-c}). \tag{14}$$

For $u\in T, l\in [\![0,|u|]\!]$ and $(k,j)\in I_K$, let $A_{u,l,k,j}$ be the number of ancestors $v\in [\![\varnothing,u[\![$ such that $d(u,v)\le l$, and for which $c_u(T)=k$ and $f_v(u)=j$.

LEMMA 8. (i) *For every $c>0$, there exists $\gamma>0$, such that, for $n$ large enough,*

$$\mathbb{P}_n(\exists(k,j)\in I_K, u\in T, |A_{u,k,j}-\mu_k|u||\ge\gamma\sqrt{|u|\log n})\le n^{-c}.$$

(ii) *For every $c>0$, there exists $\gamma>0$ such that, for $n$ large enough,*

$$\mathbb{P}_n(\exists(k,j)\in I_K, u\in T, l\in(0,|u|], |A_{u,l,k,j}-\mu_k l|\ge\gamma\sqrt{l\log n})\le n^{-c}.$$

PROOF. (ii) clearly implies (i). But let us prove (i) first. Using (9) and (13), we have for some constant $c>0$, for any $m\in[\![0,n]\!]$, any $h\ge1$, any $\mathsf{a}\in\mathbb{N}^I[h]$,

$$\mathbb{P}_n(A_{u(m)}=\mathsf{a})\le cn^{3/2}\mathbb{Q}_h(\mathsf{a})\mathbb{1}_{h\le n}. \tag{15}$$

Then

$$\mathbb{P}_n(\exists m\in[\![0,n]\!],(k,j)\in I_K, |A_{u(m),k,j}-\mu_k|u(m)||\ge\gamma\sqrt{|u(m)|\log n})$$
$$\le cn^{3/2}\sum_{m=0}^n\sum_{h=0}^n\mathbb{P}(\exists(k,j)\in I_K, |\mathcal{M}_{k,j}^{(h)}-\mu_k|h||\ge\gamma\sqrt{h\log n}).$$

This latter probability is smaller, for any $m\le n$, $h\le n$ than $\#I_K n^{-\gamma^2/2}$ by Hoeffding. Hence, $\mathbb{P}_n(\exists(k,j)\in I_K, u\in T, |A_{u,k,j}-\mu_k|u||\ge\gamma\sqrt{|u|\log n})\le cn^{7/2}n^{-2\gamma^2}$.

For (ii), assume that $u(m)=h$ and for $l\le h$, take $v_1,\dots,v_l$ the ancestors of $u(m)$ at depth $0\le h_1<\dots<h_l<h$, and set $A'_{u(m),l,k,j}=\#\{i,c_{v_i}=$ $k, f_{v_i}(u(m)) = j\}$, the lineage of $u(m)$ restricted to the nodes $v_i$'s. By "symmetry," $(A'_{u(m),l,k,j})_{k,j}$ and $(A_{u(m),l,k,j})_{k,j}$ have the same distributions. Here "symmetry" means the following: let $v_1$ and $v_2$ be two ancestors of $u(m)$. Exchange in $T$ the two nodes $v_1$ and $v_2$ together with the subtrees rooted on their children not on $[\![0, u(m)]\!]$, as on Figure 5. We get $T'$. First $T'$ and $T$ have the same weight under $\mathbb{P}_n$. Second, $A_{u(m)}$ has the same value in $T$ and $T'$, and the nodes $u(m)$ in $T$ and $T'$ have the same depth [$u(m)$ is by definition the $m$th node].

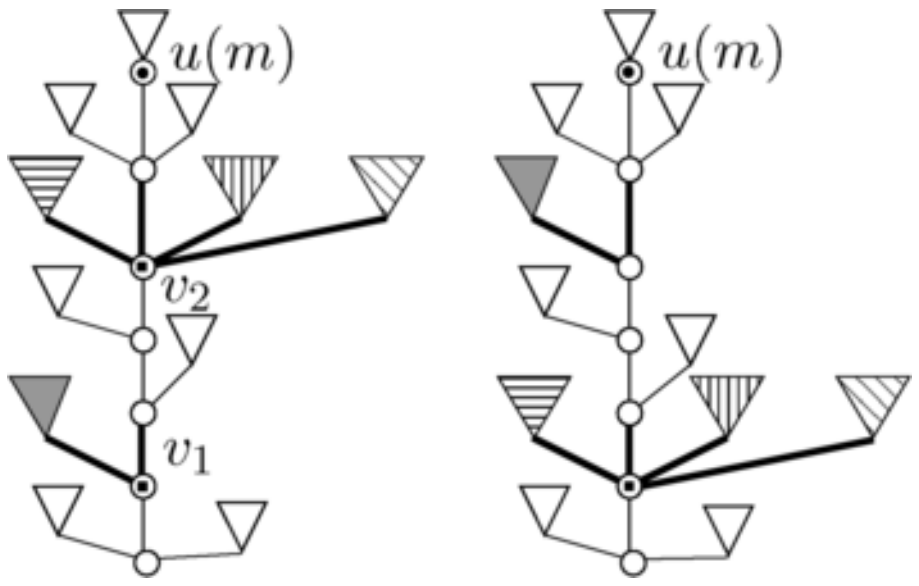

FIG. 5. *Exchange of two nodes in a lineage.*

Now take $v$ the ancestor of $u(m)$ at depth $l$. By symmetry, $(A_{v,k,j})_{k,j}$ and $(A_{u(m),l,k,j})_{k,j}$ have the same distributions. And thus, by (i), for any $m \le n$, $l \le n$, $\mathbb{P}(\exists (k,j) \in I_K, |A_{u(m),l,k,j} - \mu_k|l|| \ge \gamma\sqrt{l \log n})$ is certainly smaller than $cn^{7/2}n^{-2\gamma^2}$. As a direct consequence, $cn^{7/2+2}n^{-2\gamma^2}$ is a bound for $\mathbb{P}_n(\exists (k,j) \in I_K, u \in T, l \in (0, |u|), |A_{u,l,k,j} - \mu_k l| \ge \gamma\sqrt{l \log n})$. $\square$

We end this section with a result concerning multinomial random variables. For any $h > 0$, set

$$J_h = \left\{ \mathsf{a} \in \mathbb{N}^I[h], (N_1(\mathsf{a}), N_2(\mathsf{a})) \in \left[ \frac{\sigma_\mu^2}{2}h - h^{2/3}, \frac{\sigma_\mu^2}{2}h + h^{2/3} \right]^2 \right\}.$$

LEMMA 9. *For any $h \in \mathbb{N}$, $N_1(\mathcal{M}^{(h)})$ and $N_2(\mathcal{M}^{(h)})$ have the same law and there exists $c_1 > 0, c_2 > 0$ such that*

$$\mathbb{P}(\mathcal{M}^{(h)} \notin J_h) \le c_1 \exp(-c_2 h^{1/3}).$$

PROOF. The first assertion is easy. Writing $\{|N_1(\mathcal{M}^{(h)}) - \frac{\sigma_\mu^2 h}{2}| \ge n^{2/3}\} \subset \bigcup_{k,j}\{|\mathcal{M}_{k,j}^{(h)} - h\mu_k| \ge h^{2/3}/\#I_K\}$ [check that $\sum_{k,j}(j-1)\mu_k = \sum_{k,j}(k-j)\mu_k = \frac{\sigma_\mu^2}{2}$], by Hoeffding, one has $\mathbb{P}(|\mathcal{M}_{k,j}^{(h)} - h\mu_k| \ge h^{2/3}/\#I_K) \le 2\exp(-h^{1/3}/(2\#I_K))$. Summing this for $(k,j) \in I_K$, the result is shown to be true. $\square$

2.2.2. *A first comparison lemma.* In this section $S$ denotes a Polish space. For any r.v. $X$ taking its values in $S$, we denote by $\mathbb{P}_X$ the distribution of $X$: that is, $\mathbb{P}_X(A) = \mathbb{P}(X \in A)$ for any $A$ Borelian of $S$.

DEFINITION 1. Letting $(Y_1, Y_2, \ldots)$ and $(X_1, X_2, \ldots)$ be two sequences of r.v. taking their values in $S$ such that $\mathbb{P}_{X_n}$ is absolutely continuous with respect to $\mathbb{P}_{Y_n}$, we write $\mathbb{P}_{X_n} \prec \mathbb{P}_{Y_n}$. Let $f_n$ be a nonnegative measurable function $f_n$ such that $\mathbb{P}_{X_n} = f_n \mathbb{P}_{Y_n}$ [in other words, $\mathbb{P}_{X_n}(A) = \int_A f_n \, d\mathbb{P}_{Y_n}$ for any Borelian $A$ of $S$]; the existence of $f_n$ is ensured by the Radon–Nikodym theorem. For any $\varepsilon > 0$, let $A^n_\varepsilon := \{x, |f_n(x) - 1| < \varepsilon\}$. We say that $\mathbb{P}_{X_n}/\mathbb{P}_{Y_n} \to 1$, or $X_n // Y_n \to 1$, if for any $\varepsilon > 0$, $\mathbb{P}_{Y_n}(A^n_\varepsilon) \to 1$ (this is a convergence of $f_n$ to 1 in a weak sense).

If $X_n // Y_n \to 1$, then $\mathbb{P}_{X_n}(A^n_\varepsilon) \to 1$, and for any $B \subset A^n_\varepsilon$, $|\mathbb{P}_{Y_n}(B) - \mathbb{P}_{X_n}(B)| \leq \varepsilon \mathbb{P}_{Y_n}(B)$, therefore, $\sup_{B \text{ Borelian}} |\mathbb{P}_{X_n}(B) - \mathbb{P}_{Y_n}(B)|$, the total variation distance between $X_n$ and $Y_n$, goes to 0. Hence, the following lemma is a straightforward consequence of the Portmanteau theorem:

LEMMA 10. *If $X_n // Y_n \to 1$ and $Y_n \xrightarrow[n]{(d)} Y$, then $X_n \to_n^{(d)} Y$.*

We end this section by an argument of continuity:

LEMMA 11. *Let $(g_n)$ be a sequence of continuous functions from $S$ into a Polish space $S'$. If $X_n // Y_n \to 1$, then $g_n(X_n) // g_n(Y_n) \to 1$.*

PROOF. If $\mathbb{P}_{X_n} \prec \mathbb{P}_{Y_n}$, then so do $\mathbb{P}_{g_n(X_n)} \prec \mathbb{P}_{g_n(Y_n)}$, and then there exists a nonnegative measurable function $h_n$ such that $\mathbb{P}_{g_n(X_n)} = h_n \mathbb{P}_{g_n(Y_n)}$. As above, denote $A^n_{\varepsilon'} := \{x, |f_n(x) - 1| < \varepsilon'\}$ where $f_n$ satisfies $\mathbb{P}_{X_n} = f_n \mathbb{P}_{Y_n}$ and $\mathbb{P}_{Y_n}(A^n_{\varepsilon'}) \to 1$ and set $B^n_{\varepsilon'} = g_n(A^n_{\varepsilon'})$. We have $\mathbb{P}_{g_n(Y_n)}(B^n_{\varepsilon'}) \to 1$. For any $A \subset B^n_{\varepsilon'}$,

$$|\mathbb{P}_{g_n(X_n)}(A) - \mathbb{P}_{g_n(Y_n)}(A)| = \left| \int_A (h_n - 1) \, d\mathbb{P}_{g_n(Y_n)} \right| < \varepsilon',$$

the inequality follows that $g_n^{-1}(A) \subset A^n_{\varepsilon'}$. Letting now $\varepsilon > 0$ be fixed and setting $A = \{x, h_n(x) - 1 > \varepsilon\} \cap B^n_{\varepsilon'}$ (or $A = \{x, h_n(x) - 1 < -\varepsilon\} \cap B^n_{\varepsilon'}$), we get $\mathbb{P}_{g_n(Y_n)}(\{x, |h_n(x) - 1| \geq \varepsilon\}) \leq 2\varepsilon'/\varepsilon$; choosing $\varepsilon' > 0$ small, one sees that this is arbitrarily small for $n$ large enough. □

2.2.3. *Proof of the convergence of the uni-dimensional distributions in Theorem 2.* In this section we work under $\mathbb{P}_n$. Let $X^n_m := (A_{u(m)}, |u(m)|)$ and $Y^n_m := (A^\star_m, |u(m)|)$, where the distribution of $A^\star_m$ knowing $|u(m)| = h$ is simply $\mathbb{Q}_h$. The aim of this section is to compare $X^n_m$ with $Y^n_m$ and to deduce from the asymptotic behavior of $Y_m^n$ the one of $X_m^n$. The proof of the convergence of the finite-dimensional distributions will also use this strategy.

For $M>0$, and $n\in\mathbb{N}$, consider

$$\Lambda_{n,M}=\{(\mathsf{a},h),h\in\sqrt{n}[M^{-1},M],\mathsf{a}\in J_h\}.$$

We have the following:

PROPOSITION 12. (i) *For any $m,n$, with $m\le n$, we have $\mathbb{P}_{X_m^n}\prec\mathbb{P}_{Y_m^n}$.*

(ii) *For any $s\in(0,1)$, $\alpha>0$, there exists $M_0$ s.t. for $n$ large enough, $\mathbb{P}_n(Y_{\lfloor ns\rfloor}^n\in\Lambda_{n,M_0})\ge 1-\alpha$ and for any $M>0$,*

$$\sup_{(a,h)\in\Lambda_{n,M}}\left|\frac{\mathbb{P}_n(X_{\lfloor ns\rfloor}^n=(\mathsf{a},h))}{\mathbb{P}_n(Y_{\lfloor ns\rfloor}^n=(\mathsf{a},h))}-1\right|\underset{n}{\to}0. \tag{16}$$

(iii) *For any $s\in(0,1)$, $X_{\lfloor ns\rfloor}^n//Y_{\lfloor ns\rfloor}^n\to 1$.*

PROOF. (iii) is a consequence of (ii). Let $\mathsf{a}\in\mathbb{N}_I[h]$. Since $\{A_{u(m)}=\mathsf{a}\}\subset\{|u(m)|=h\}$, $\mathbb{P}_n((A_{u(m)},|u(m)|)=(\mathsf{a},h))=\mathbb{P}_n(A_{u(m)}=\mathsf{a})$. According to Proposition 5 and formula (10),

$$\frac{\mathbb{P}_n(X_m^n=(\mathsf{a},h))}{\mathbb{P}_n(Y_m^n=(\mathsf{a},h))}=\frac{\mathbb{P}(|\mathbf{f}_{N_1(\mathsf{a})}|=m-h,|\mathbf{f}'_{1+N_2(\mathsf{a})}|=n+1-m)}{\mathbb{P}(|\mathbf{f}_{N_1(\mathcal{M}^{(h)})}|=m-h,|\mathbf{f}'_{1+N_2(\mathcal{M}^{(h)})}|=n+1-m)}. \tag{17}$$

Then (i) holds true. Assume now that $s\in(0,1)$ and $\alpha>0$ are fixed. There exists $M$ such that, for $n$ large enough, $\mathbb{P}_n(|u(\lfloor ns\rfloor)|\in\sqrt{n}[M^{-1},M])\ge 1-\alpha/2$ [since $\mathbf{h}_n\overset{(d)}{\underset{n}{\to}}\frac{2}{\sigma_\mu}\mathbf{e}$ and since $\mathbb{P}(\mathbf{e}_s=0)=0$ for any $s\in(0,1)$]. For such a $M$,

$$\begin{aligned}
&\mathbb{P}_n(Y_{\lfloor ns\rfloor}^n\in\Lambda_{n,M})\\
&\quad=\mathbb{P}_n(Y_{\lfloor ns\rfloor}^n\in\Lambda_{n,M},|u(\lfloor ns\rfloor)|\in\sqrt{n}[M^{-1},M])\\
&\quad=\sum_{l\in\sqrt{n}[M^{-1},M]}\mathbb{P}_n(|u(\lfloor ns\rfloor)|=l)\mathbb{P}_n(Y_{\lfloor ns\rfloor}^n\in\Lambda_{n,M}||u(\lfloor ns\rfloor)|=l)\\
&\quad\ge\mathbb{P}_n(|u(\lfloor ns\rfloor)|\in\sqrt{n}[M^{-1},M])\min_{l\in\sqrt{n}[M^{-1},M]}\mathbb{P}(\mathcal{M}^{(l)}\in J_l).
\end{aligned}$$

This minimum goes to 1, thanks to Lemma 9.

Now, according to Lemma 6(i) and (ii), since $\mathbf{f}$ and $\mathbf{f}'$ are independent,

$$\begin{aligned}
&\mathbb{P}(|\mathbf{f}_{N_1(\mathsf{a})}|=\lfloor ns\rfloor-h,|\mathbf{f}'_{1+N_2(\mathsf{a})}|=n+1-\lfloor ns\rfloor)\\
&\quad=\frac{N_1(\mathsf{a})(1+N_2(\mathsf{a}))}{(\lfloor ns\rfloor-h)(n+1-\lfloor ns\rfloor)}\\
&\qquad\times\mathbb{P}(W_{\lfloor ns\rfloor-h}=-N_1(\mathsf{a}))\mathbb{P}(W_{n-\lfloor ns\rfloor+1}=-N_2(\mathsf{a})-1),
\end{aligned}$$

and then for any $M>0$,

$$\sup_{(a,h)\in\Lambda_{n,M}}\left|\frac{\mathbb{P}(|\mathbf{f}_{N_1(\mathsf{a})}|=\lfloor ns\rfloor-h,|\mathbf{f}'_{1+N_2(\mathsf{a})}|=n+1-\lfloor ns\rfloor)}{q_{n,s,h}}-1\right|\underset{n}{\rightarrow}0$$

for

$$q_{n,s,h}=\frac{\sigma_\mu^2h^2\exp(-\sigma_\mu^4h^2/(8ns(1-s)))}{8\pi n^3(s(1-s))^{3/2}}.$$

Now, $\mathbb{P}(|\mathbf{f}_{N_1(\mathcal{M}^{(h)})}|=\lfloor ns\rfloor-h,|\mathbf{f}'_{1+N_2(\mathcal{M}^{(h)})}|=1+n-\lfloor ns\rfloor)=A_h+B_h$, where

$$A_h:=\mathbb{P}(|\mathbf{f}_{N_1(\mathcal{M}^{(h)})}|=\lfloor ns\rfloor-h,|\mathbf{f}'_{1+N_2(\mathcal{M}^{(h)})}|=1+n-\lfloor ns\rfloor,\mathcal{M}^{(h)}\notin J_h)$$

$$B_h:=\mathbb{P}(|\mathbf{f}_{N_1(\mathcal{M}^{(h)})}|=\lfloor ns\rfloor-h,|\mathbf{f}'_{1+N_2(\mathcal{M}^{(h)})}|=1+n-\lfloor ns\rfloor,\mathcal{M}^{(h)}\in J_h).$$

Using again Lemma 6(i) and (ii), we get

$$\sup_{h\in\sqrt{n}[M^{-1},M]}\left|\frac{B_h}{q_{n,s,h}}-1\right|\underset{n}{\rightarrow}0.$$

On the other hand, $A_h\le\mathbb{P}(\mathcal{M}^{(h)}\notin J_h)\le c_1\exp(-c_2h^{1/3})\le 2\exp(-cn^{1/6}/M)$ for any $h\in\sqrt{n}[M^{-1},M]$. To complete the proof of (ii), check that $\sup_{h\in\sqrt{n}[M^{-1},M]}|A_h/B_h|\underset{n}{\rightarrow}0$. □

We have now all the tools to conclude the following:

Corollary 13. *For any $s\in(0,1)$, letting $s_n=\lfloor ns\rfloor/n$, we have*

$$(\mathbf{G}^{(n)}(s_n),\mathbf{h}_n(s_n))//(\mathcal{G}(n,\sqrt{n}\mathbf{h}_n(s_n)),\mathbf{h}_n(s_n))\underset{n}{\rightarrow}1,$$

*and the convergence of the uni-dimensional distributions holds in Theorem 2.*

(Recall that $\mathcal{G}$ is defined in Section 2.1.)

Proof of Corollary 13. Proposition 12 and Lemma 11 yield the first assertion of the Corollary.

For the second assertion, we first examine $s=0$ and $s=1$. Since $\mathbf{h}_n(0)=0$ and Remark 1 entails that $\mathbf{h}_n(1)\overset{\text{proba.}}{\underset{n}{\rightarrow}}0$, the convergence of the uni-dimensional distributions holds in Theorem 2 for $s=0$ and $s=1$.

For $s\in(0,1)$, since $\mathbf{h}_n\overset{(d)}{\underset{n}{\rightarrow}}\mathbf{h}$ in $C[0,1]$, by the Skorohod representation theorem [14], Theorem 3.30, there exists a probability space on which this convergence is a.s. On this space [or on an augmented space on which the pair $(\mathcal{G}(n,\sqrt{n}\mathbf{h}_n(s_n)),\mathbf{h}_n(s_n))$] is defined

$$(\mathcal{G}(n,\sqrt{n}\mathbf{h}_n(s_n)),\mathbf{h}_n(s_n))\overset{(d)}{\underset{n}{\rightarrow}}(\mathcal{G}_\infty^{\mathbf{h}_s},\mathbf{h}_s),\tag{18}$$

where $\mathcal{G}_\infty^{\mathbf{h}_s}$ is a Gaussian process which covariance function [see (8)] allows to check that $(\mathcal{G}_\infty^{\mathbf{h}_s}, \mathbf{h}_s) \stackrel{(d)}{=} (\mathbf{G}(s), \mathbf{h}_s)$ for any $s \in (0,1)$. To prove that the convergence of the uni-dimensional distribution holds in Theorem 2, it remains to control the distance between $(\mathbf{G}^{(n)}(s_n), \mathbf{h}_n(s_n))$ and $(\mathbf{G}^{(n)}(s), \mathbf{h}_n(s))$. Since $\mathbf{h}_n \xrightarrow[n]{(d)} \mathbf{h}$ in $C[0,1]$, $|\mathbf{h}_n(s_n) - \mathbf{h}_n(s)| \xrightarrow[n]{\text{proba.}} 0$. For $\mathbf{G}^{(n)}$, this is more complex, and we will establish some bounds useful also for the proof of the tightness. Let

$$\Omega_n^\rho = \Big\{ T \in \mathcal{T}_n, \max_l ||u(l+1)| - |u(l)|| \le \rho \log n \Big\}.$$

Let $\varepsilon > 0$. According to Lemma 7, for $\rho$ large enough, $\mathbb{P}_n(\Omega_n^\rho) > 1 - \varepsilon$ for $n$ large enough. We have, for $s'_n = \lfloor ns + 1 \rfloor / n$,

$$\|\mathbf{G}^{(n)}(s_n) - \mathbf{G}^{(n)}(s)\|_1 \mathbb{1}_{\Omega_n^\rho} = n(s - s_n) \sum_{i \in I_K} \mathbb{1}_{\Omega_n^\rho} |\mathbf{G}_i^{(n)}(s'_n) - \mathbf{G}_i^{(n)}(s_n)|. \tag{19}$$

In $\Omega_n^\rho$, for any $k, j$, the differences $|\mathbf{G}_{k,j}^{(n)}(s'_n) - \mathbf{G}_{k,j}^{(n)}(s_n)|$ are bounded by $2\rho n^{-1/4} \log n$ (which is a bound on the number of noncommon ancestors of two consecutive nodes in the LO for a tree in $\Omega_n^\rho$). Hence, since $s - s_n \le 1/n$, for any $\varepsilon' > 0$, for $n$ large enough,

$$\|\mathbf{G}^{(n)}(s_n) - \mathbf{G}^{(n)}(s)\|_1 \mathbb{1}_{\Omega_n^\rho} \le c(s - s_n)^{1/4 - \varepsilon'} \tag{20}$$

for some constant $c$. One concludes that $\|\mathbf{G}^{(n)}(s_n) - \mathbf{G}^{(n)}(s)\|_1 \xrightarrow[n]{\text{proba.}} 0$. □

2.3. *Convergence of the finite-dimensional distributions.* In this Section $\kappa \ge 2$ is a fixed integer. We denote by $\mathbf{s}^{(\kappa)}$ the vector $(s_1, \ldots, s_\kappa)$ where $0 < s_1 < \cdots < s_\kappa \le 1$ are fixed. Let $T \in \mathcal{T}_n$. For $i \in [\![1, \kappa]\!]$, set $u_i = u(\lfloor ns_i \rfloor)$, $u_0 = u_{\kappa+1} = \varnothing$, and

$$L(T) = \{u_i, i \in [\![1, \kappa]\!]\}.$$

We assume that $n$ is large enough such that $\lfloor ns_1 \rfloor \ge 1$, and $\lfloor ns_i \rfloor \ne \lfloor ns_j \rfloor$ for $i \ne j$, so that the $u_i$'s are different nodes of $T$ sorted according to the LO.

The aim of this section is to study the distribution of $(A_{u_i})_{i \in [\![1,\kappa]\!]}$ under $\mathbb{P}_n$, and to deduce from this the convergence of the finite-dimensional distribution in Theorem 2. The ideas are of the same type as in the case of the uni-dimensional distributions, but the details are more involved since the dependences between the r.v. $A_{u_i}$'s must be taken into account. For this, the shape of the tree spanned by the $u_i$'s must be considered.

Denote by $\check{u}_{i,j}$ the deepest (i.e., youngest) common ancestor of $u_i$ and $u_j$. Let $T_{\mathbf{s}^{(\kappa)}} = \bigcup_{i=1}^{\kappa} [\![\varnothing, u_i]\!]$ be the subtree "spanned" by the $u_i$'s,

$$Z(T) = \{\check{u}_{i,j}, 1 \le i < j \le \kappa\} = \{\check{u}_{i,i+1}, i \in [\![1, \kappa - 1]\!]\},$$

$$Z^\star(T) = Z(T) \cap L(T),$$

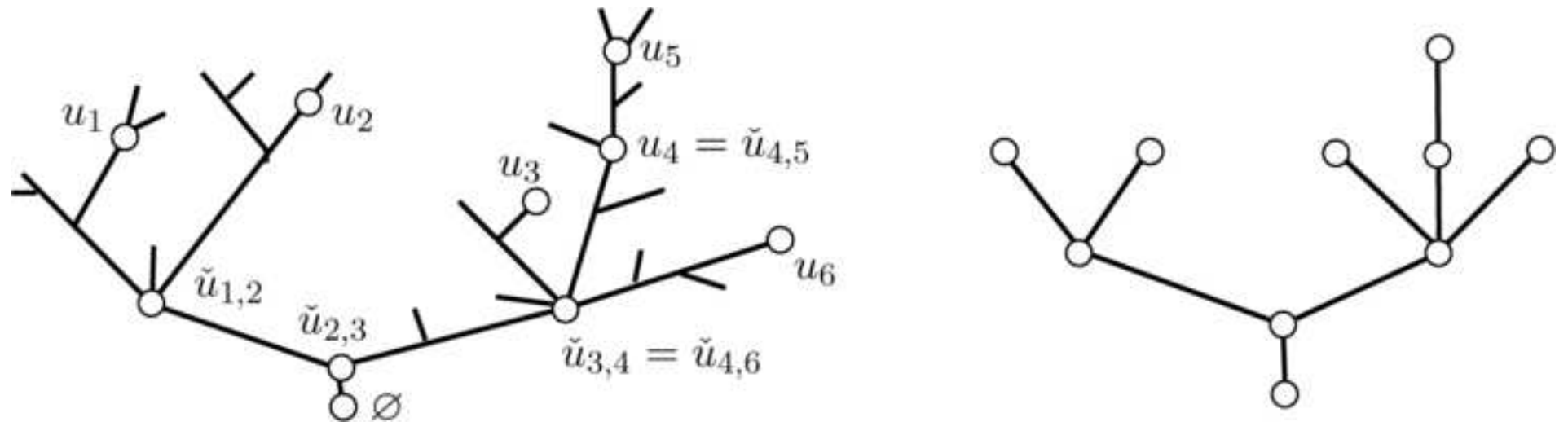

FIG. 6. *A tree $T$ and the associated tree $\Phi(T) = \{\varnothing, 1, 11, 111, 112, 12, 121, 122, 1221, 123\}$. On this example $Z^\star(T) = u_4$, and then $I_T^{b(n)} = \{122\}$.*

the set of *branching nodes* in $T_{\mathbf{s}^{(\kappa)}}$, and the nodes in $L(T)$ that are ancestors of other nodes of $L(T)$.

DEFINITION 2. The shape function $b(n)$ associates with $T_{\mathbf{s}^{(\kappa)}}$ the smallest tree in $\mathcal{T}$ having the same branching structure (that we call shape) together with a coding of the nodes of $Z^\star(T)$ (see Figure 6). Formally,

$$b(n) : \mathcal{T}_n \times [0,1]^\kappa \longrightarrow \mathcal{T} \times \mathcal{P}_F(\mathbb{U}),$$
$$(T, \mathbf{s}^{(\kappa)}) \longmapsto (T^{b(n)}, I_T^{b(n)}),$$

where $\mathcal{P}_F(\mathbb{U})$ is the set of finite subsets of $\mathbb{U}$ and where $T^{b(n)}$ is characterized by:

(i) $T^{b(n)}$ is a tree having $\#(Z(T) \cup L(T) \cup \{\varnothing\})$ nodes,

(ii) there exists an increasing function $\Phi_T$ from $Z(T) \cup L(T) \cup \{\varnothing\}$ in $T^{b(n)}$, preserving the descendants: $\Phi_T(u)$ is an ancestor of $\Phi_T(v)$ in $T^{b(n)}$ iff $u$ is an ancestor of $v$ in $T$.

The set $I_T^{b(n)}$ is defined to be $\Phi_T(Z^\star(T))$.

The tree $T^{b(n)}$ can be constructed in somehow squeezing the paths between the nodes of $Z(T) \cup L(T) \cup \{\varnothing\}$ in unit length edge and in renaming the vertices in order to get a tree. The function $\Phi_T$ is unique and for short, for any $u \in Z(T) \cup L(T) \cup \{\varnothing\}$, we write $u^{b(n)}$ instead of $\Phi_T(u)$. The set $I_T^{b(n)}$ encodes the images of the nodes of $Z^\star(T)$. Notice that $\# I_T^{b(n)} = \kappa - \# \partial T^{b(n)}$, and when $I_T^{b(n)}$ is not empty, the tree $T^{b(n)}$ alone is not sufficient in general to guess $I_T^{b(n)}$.

In what follows, we will often write $b$ instead of $b(n)$.

A pair $(u, v)$ [with $u, v \in L(T) \cup Z(T)\{\varnothing\}$] such that $u^b = fa(v^b)$ is the father of $v^b$ in $T^b$ will be called a spanned branch. The contents of the spanned branches will be carefully handled since they contribute in general to several $A_{u_i}$'s. The set of spanned branches can naturally be indexed by the edges $(u^b, v^b)$ of $T^b$, but also by the nodes $v^b$ of $T^b \setminus \{\varnothing\}$ using the bijection between the edges of $T^b$ and $T^b \setminus \{\varnothing\}$ that associates with the edge $(u^b, v^b)$ the node $v^b$.

Using this labeling, we define $A_{(v^b)}$ the content of the spanned edge $(u, v)$ by

$$A_{(v^b),k,j} := \#\{w \in ]\!]u, v[\![, c_w = k, f_w(v) = j\}.$$

The extremities of the spanned branches are not counted in the $A_{(v^b),k,j}$'s in order to simplify the decompositions ($\varnothing$ was counted in the unidimensional case). It is easy to check that

$$A_{(v^b),k,j} = (A_{v,k,j} - A_{u,k,j}) - \mathbb{1}_{(c_u, f_u(v))}(k, j). \tag{21}$$

We also introduce the "ordered content" of the edges. For any $v^b \in T^b \setminus \{\varnothing\}$, define $\overrightarrow{A}_{(v^b)}(T)$, the ordered content of the edge $(u, v)$ by

$$\overrightarrow{A}_{(v^b)}(T) := ((c_w, f_w(v)), w \in ]\!]u, v[\![),$$

the nodes of $]\!]u, v[\![$ being sorted according to the LO.

We write simply $\overrightarrow{A}(T) = (\overrightarrow{A}_{(v^b)}(T))_{v^b \in T^b \setminus \{\varnothing\}}$, the list of ordered contents of the spanned edges.

The ordered content of any edge belongs to $\bigcup_{i \geq 0} (I_K)^i$. The canonical surjection $\pi$ from $\bigcup_{i \geq 0} (I_K)^i$ into $\mathbb{N}^{I_K}$ associates with the ordered content $\overrightarrow{B} = ((k_i, j_i), i = 1, \ldots, l)$ the content $B$:

$$B_{k,j} = (\pi(\overrightarrow{B}))_{k,j} = \#\{i, (k_i, j_i) = (k, j)\}.$$

The application $\pi$ can be extended to the list of ordered contents, and we set

$$A(T) = \pi((\overrightarrow{A}_{(v^b)}(T)))_{v^b \in T^b \setminus \{\varnothing\}} = (\pi(\overrightarrow{A}_{(v^b)}(T)))_{v^b \in T^b \setminus \{\varnothing\}}.$$

The definition of $N_1$ and $N_2$ (given in Section 2.2) are extended to $\bigcup_{i \geq 0} (I_K)^i$: we set,

$$N_1(\overrightarrow{B}) := N_1(\pi(\overrightarrow{B})) \quad \text{and} \quad N_2(\overrightarrow{B}) := N_2(\pi(\overrightarrow{B})).$$

We denote by $H_{v^b}$ the cardinality $\#]\!]u, v[\![$ where $u^b = fa(v^b)$. We set

$$\mathcal{H}_T = (H_{v^b})_{v^b \in T^b \setminus \{\varnothing\}},$$

the ordered list of the spanned branches lengths.

The nodes of $Z(T) \cup \{\varnothing\}$ are the "hinge nodes" laying between the spanned edges, and they also contribute to the $A_{u_i}$'s. For any $u$ in $Z(T) \cup \{\varnothing\}$, the set $f_u(L(T))$ is a subset of $[\![1, \ldots, c_u(T)]\!]$ with $c_{u^b}(T^b)$ elements: the set of the ranks of the children of $u$ that are ancestors of the nodes of $L(T)$. We encode the contribution of $Z(T) \cup \{\varnothing\}$ to the lineage, thanks to the sequence $\Theta_T$:

$$\Theta_T = (C(u^b), R(u^b 1), \ldots, R(u^b C(u^b)))_{u^b \in T^b \setminus \partial T^b},$$

where $C(u^b) = c_u(T)$ and $R(u^b 1) < \cdots < R(u^b C(u^b))$ is the sorted list of the elements of $f_u(L(T))$. Note that $u^b 1, \ldots, u^b C(u^b)$ are the children of $u^b$ in $T^b$ and then the arguments of $R$ are unambiguous.

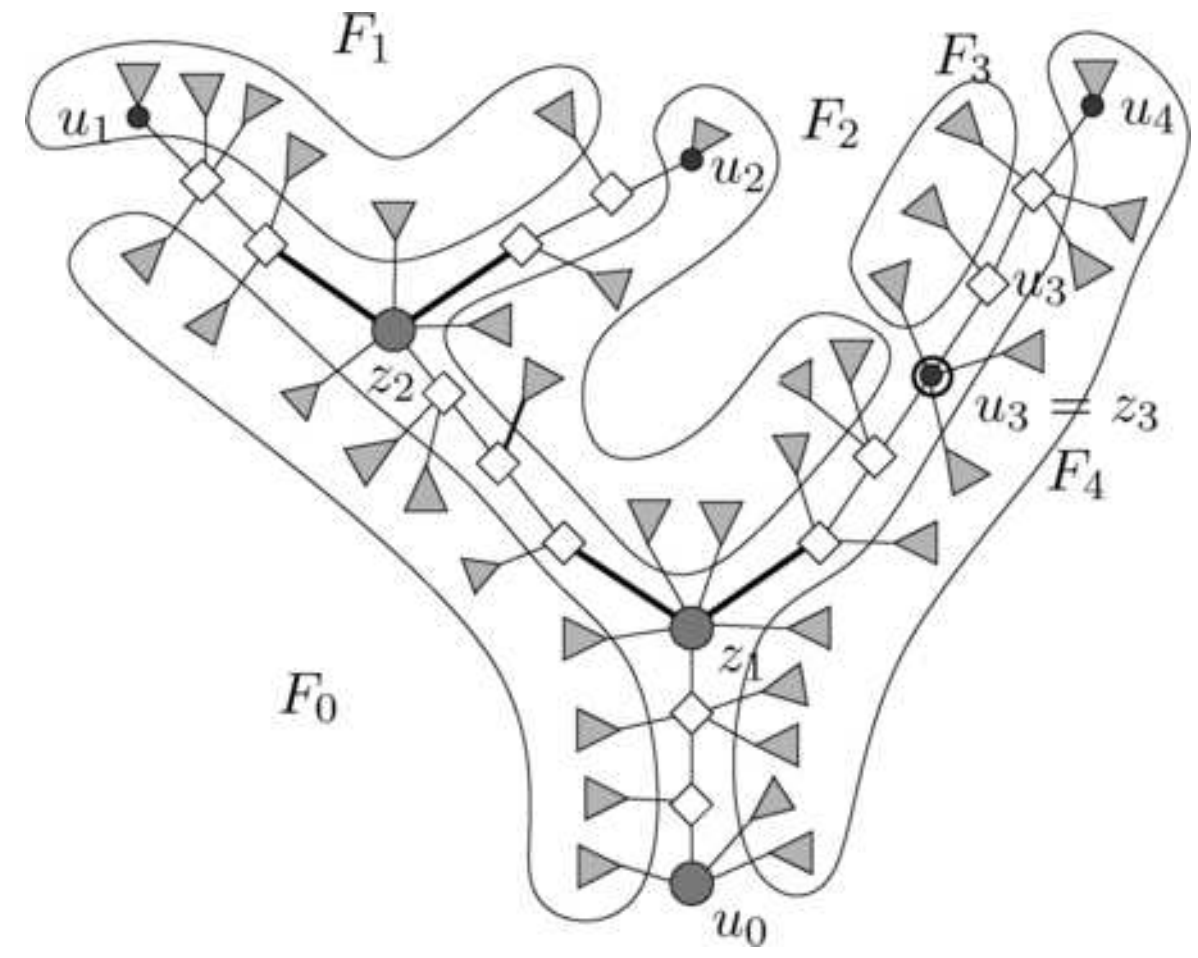

FIG. 7. *The forests considered in the decomposition are surrendered. Notice that the node $u_i$ belong to a forest only if it is a leaf in $T_{s(\kappa)}$. Observe also the contribution of the neighbors of the "hinge nodes," the nodes $z_i$'s on the picture.*

The idea now is the following. If $(T^{b(n)}, I_T^{b(n)}, \overrightarrow{A}(T), \Theta_T)$ is known, to end the description of $T$ using $T_{s(\kappa)}$, it remains to describe the fringe subtrees rooted in the neighborhood of $T_{s(\kappa)}$ (the *fringe subtree* of $T$ rooted at $u$ is $T_u = \{v \in \mathbb{U} : uv \in T\}$). We pack these subtrees into forests that are, up to some border effects, rooted on the neighbors of $T_{s(\kappa)}$ between $u_i$ and $u_{i+1}$.

For any simple path $I$ in $T$, we denote by $\mathcal{N}(I)$ the neighborhood of $I$:

$$\mathcal{N}(I) := \{u \in T, d_T(u, I) = 1\}.$$

We now build the set of roots of the forest we consider (see Figure 7):

$$\mathcal{S}_0 = \{v \in \mathcal{N}(\llbracket \varnothing, u_1 \rrbracket), \varnothing \prec v \prec u_1\},$$

$$\mathcal{S}_i = \begin{cases} \{v \in \mathcal{N}(\rrbracket u_i, u_{i+1} \rrbracket), u_i \prec v \prec u_{i+1}\} \cup \{u_i\}, & \text{if } u_i \in L(T) \setminus Z^\star(T), \\ \{v \in \mathcal{N}(\llbracket u_i, u_{i+1} \rrbracket), u_i \prec v \prec u_{i+1}\}, & \text{if } u_i \in Z^\star(T), \end{cases}$$

$$i \in \llbracket 1, \kappa - 1 \rrbracket,$$

$$\mathcal{S}_\kappa = \{v \in \mathcal{N}(\rrbracket u_\kappa, \varnothing \rrbracket), u_\kappa \prec v\} \cup \{u_\kappa\}.$$

The forests we consider are

$$\mathfrak{F}_i(T) := (T_u, u \in \mathcal{S}_i);$$

we denote by $\mathfrak{F}(T) = (\mathfrak{F}_i(T))_{i=0,\dots,\kappa}$ the (ordered) sequence of forests.

Let $l \in [\![0, \kappa-1]\!]$ be fixed. Let $u_l^b = v_0^b, v_1^b, \dots, v_m^b = u_{l+1}^b$, the shortest path in $T^b$ between $u_l^b$ and $u_{l+1}^b$. Let $v_i^b = u_l^b \vee u_{l+1}^b$ be the deepest (youngest) common ancestor of $u_l^b$ and $u_{l+1}^b$. Since $v_0^b \prec v_m^b$, two cases arise, $v_i^b = v_0^b$ or $0 < i < m$. We have

$$\begin{aligned} &\#\mathcal{S}_l(T^b, I_T^b, \overrightarrow{A}_T, \Theta(T)) \\ &\qquad = \mathcal{N}_{l,l+1}(T^b, I_T^b, \overrightarrow{A}(T)) + \mathcal{Y}_{l,l+1}(T^b, I_T^b, \Theta(T)), \end{aligned} \tag{22}$$

where, for any $l$, $\mathcal{N}_{l,l+1}(T^b, I_T^b, \overrightarrow{A}(T))$ counts the number of subtrees rooted on the neighbors of the spanned branches visited between $u_l$ and $u_{l+1}$, according that these subtrees are on the right or on the left of these spanned branches:

$$\mathcal{N}_{l,l+1}(T^b, I_T^b, \overrightarrow{A}(T)) = \sum_{p=0}^{i-1} N_2(\overrightarrow{A}_{(v_p)}) + \sum_{p=i+1}^{m} N_1(\overrightarrow{A}_{(v_p)}),$$

and $\mathcal{Y}_{l,l+1}(T^b, I_T^b, \Theta_T)$ counts the number of subtrees rooted on the neighbors of the nodes of $Z(T) \cup \{\varnothing\}$:

$$\begin{aligned} \mathcal{Y}_{l,l+1}(T^b, I_T^b, \Theta_T) &= \mathbb{1}_{u_l^b \in \partial T^b} + R(v_{i+1}^p) - R(v_{i-1}^p)\mathbb{1}_{i>0} - 1 \\ &\quad + \sum_{p=1}^{i-1}[C(v_p^b) - R(v_p^b C(v_p^b))] + \sum_{p=i+1}^{m-1}[R(v_p^b) - 1]. \end{aligned}$$

Indeed, $R(v_{i+1}^p) - R(v_{i-1}^p)\mathbb{1}_{i>0} - 1$ is the number of children of $v_i^p$ in $\mathcal{S}_l$, the sum $\sum_{p=1}^{i-1}[C(v_p^b) - R(v_p^b C(v_p^b))]$ counts the number of children of $v_1^p, \dots, v_{i-1}^p$ in $S_l$, and $\sum_{p=i+1}^{m-1}[R(v_p^b) - 1]$ the children of $v_{i+1}^p, \dots, v_{m-1}^p$ in $S_l$.

For $l = \kappa$, let $u_\kappa^b = v_0^b, v_1^b, \dots, v_m^b = u_{\kappa+1}^b = \varnothing$. Since $v_i^b$ is an ancestor of $v_{i-1}^b$, we have

$$\mathcal{N}_{l,l+1}(T^b, I_T^b, \overrightarrow{A}(T)) = \sum_{p=0}^{m-1} N_2(\overrightarrow{A}_{(v_p)}),$$

$$\mathcal{Y}_{l,l+1}(T^b, I_T^b, \Theta_T) = \sum_{p=1}^{m}[C(v_p^b) - R(v_p^b C(v_p^b))],$$

notice that the term $p = m$ concerns the root.

The cardinalities $F_l(T) = \#\mathfrak{F}_i(T)$ satisfy

$$F_l(\mathcal{H}_T, T^b) = (\lfloor ns_{l+1} \rfloor - \lfloor ns_l \rfloor + 1) - (|u_{l+1}| - |\check{u}_{l,l+1}|) - \mathbb{1}_{u_l^b \notin \partial T^b}, \tag{23}$$

since the visit times of the nodes $u_i$ are $\lfloor ns_i \rfloor$ and since $|u_{l+1}| - |\check{u}_{l,l+1}| + 1$ nodes visited during $[\![ns_l, ns_{l+1}]\!]$ are not in $\bigcup_{u \in \mathcal{S}_l} T_u$.

For any tree $t$ having $n$ nodes, $\mathbb{P}(\mathbf{T} = t) = \prod_{u \in t} p_{c_u(t)}$. Hence, if $t$ has $n$ nodes, putting together the contribution of the forests and of the ordered edge contents, we get

$$\begin{aligned} \mathbb{P}(\mathbf{T} = t) = &\left( \prod_{i=0}^{l} \mathbb{P}(\mathbf{f}_{\#\mathcal{S}(t^b, A(t), \Theta(t))} = \mathfrak{F}_i(t)) \right) \\ &\times \left( \prod_{v^b \in t^b \setminus \{\varnothing\}} \prod_{k,j} p_k^{\pi(\overrightarrow{A}_{(v^b)}(t))_{k,j}} \right) \left( \prod_{z^b \in t^b \setminus \partial t^b} p_{C(z^b)} \right), \end{aligned} \tag{24}$$

where $\mathbf{f}_x$ is a random forest with $x$ roots (see Section 2.2).

Comment 1. A few points worth mentioning:

– Let $A = (A_{k,j})_{(k,j) \in I_K}$ with $\sum_{i \in I_K} A_i = h$. The number of ordered contents having $A$ as content is

$$\#\pi^{-1}(A) = \binom{h}{(A_i)_{i \in I_K}}, \tag{25}$$

and thus,

$$\sum_{\overrightarrow{B} \in \pi^{-1}(A)} \prod_{k,j} p_k^{\pi(\overrightarrow{B})_{k,j}} = \binom{h}{(A_i)_{i \in I_K}} \prod_{k,j} p_k^{A_{k,j}} = \mathbb{Q}_h(A). \tag{26}$$

This is simply due to the fact that all permutations of the "symbols" $(k, j)$'s are possible in the ordered contents $\overrightarrow{B}$'s such that $\pi(\overrightarrow{B}) = A$.

– The size of the forests $\mathfrak{F}_i$, as well as their number of trees, is a function of the contents (it does not depend on the order of the content).

We are now able to express the probability to observe a shape together with the (ordered) contents in terms of the probability that some forests have some prescribed sizes.

The probability to observe some contents will be obtained by summing on all corresponding ordered contents, thanks to (26). We won't do this job on all possible shapes since asymptotically only "the simplest shapes" eventually happen. We first state a result in this direction.

2.3.1. *Decomposition of a tree $T$ given $\mathcal{A}_T$ and $T^b$.* Let $\mathcal{T}^B_{2\kappa-1} = \{T \in \mathcal{T}_{2\kappa-1}, \deg(\varnothing) = 1, \forall u \in T \setminus \varnothing, \deg(u) \in \{0, 2\}\}$ be the set of trees with $2\kappa - 1$ edges, with binary branching points (except the root that has only one child). Denote by

$$\begin{aligned} \Delta_{n,M} = \{&T \in \mathcal{T}_n, T^b \in \mathcal{T}^B_{2\kappa-1}, \forall (u^b, v^b) \in T^b, \\ &u^b = fa(v^b), d(u,v) \in \sqrt{n}[M^{-1}, M]\}. \end{aligned}$$

A tree in $\Delta_{n,M}$ has its shape in $\mathcal{T}^B_{2\kappa-1}$, and all its spanned branches lengths in $\sqrt{n}[M^{-1}, M]$. By a counting argument, $I^b_T$ is empty, which means that the nodes $u_i$ are sent on the leaves of $T^b$ by $\Phi_T$.

LEMMA 14. *For any $\varepsilon > 0$, there exists $M > 0$ such that, for $n$ large enough,*

$$\mathbb{P}_n(\Delta_{n,M}) \geq 1 - \varepsilon.$$

PROOF. This is a consequence of $\mathbf{h}_n \xrightarrow[n]{(d)} \mathbf{h} = 2\mathsf{e}/\sigma_\mu$ and of the properties of $\mathsf{e}$: $\mathsf{e}$ is a.s. nonnull on $(0,1)$ and the local minima of $\mathsf{e}$ are a.s. all different (the continuum random tree is a.s. a binary tree). $\square$

For any $T \in \Delta_{n,M}$, $\Phi_T$ sends the nodes of $L(T)$ on the leaves of $T^b$ and the nodes of $Z(T)$ on the internal nodes of $T^b \setminus \{\varnothing\}$. The branching points are distinct $\#Z(T) = \kappa - 1$ and $L(T) \cap Z(T)$ is empty. Therefore, $\Theta_T$ belongs to $\mathcal{D}^b = \mathcal{D}_2 \times \mathcal{D}_3^{\kappa-1}$, where $\mathcal{D}_2 = \{(c,x), 1 \leq x \leq c\}$ and $\mathcal{D}_3 = \{(c,x,y), 1 \leq x < y \leq c\}$. Note that, under (H$_1$), $\mathcal{D}^b$ is a subset of $[\![1, K]\!]^{3\kappa-1}$.

2.3.2. *A second comparison result.* The idea is to compare $(A(\mathbf{T}), \mathcal{H}_{\mathbf{T}}, \Theta_{\mathbf{T}}, \mathbf{T}^b)$ under $\mathbb{P}_n$ with some simple random variables. Thanks to the previous lemma, we will somehow consider only the case $I^b_T$ is empty.

Denote by $(\mathcal{A}_n, \mathcal{H}_n, \Theta_n, \mathbf{T}^b_n)$ the r.v. $(A(\mathbf{T}), \mathcal{H}_{\mathbf{T}}, \Theta_{\mathbf{T}}, \mathbf{T}^b)$ under $\mathbb{P}_n$. Hence, $\mathcal{A}_n = (\mathcal{A}^i_n)_{i=1,\dots,K}$ is the sequence of contents, $\mathcal{H}_n = (\mathcal{H}^1_n, \dots, \mathcal{H}^K_n)$ the sequence of spanned edge lengths, $\Theta_n$ the sequence of branching properties of the hinge nodes where $K$ is then the number of spanned branches (the number of edges of $\mathbf{T}^b_n$). These sequences were labeled by the nodes of the tree $\mathbf{T}^b$, but, we may and will consider that they are labeled by integers (the LO is a total order). This is equivalent knowing the shape and allows one to work also when the shape is not known.

We define now the 4-tuple $(\mathcal{A}^\star_n, \mathcal{H}_n, \Theta^\star_n, \mathbf{T}^b_n)$ as follows: $\mathcal{H}_n$ and $\mathbf{T}^b_n$ have the same law as above, and $\mathcal{A}^\star_n$ and $\Theta^\star_n$ are described conditionally on $\mathcal{H}_n$ and $\mathbf{T}^b_n$. Conditionally on $\mathcal{H}_n = (\mathcal{H}^1_n, \dots, \mathcal{H}^K_n)$, where $K$ is then the number of spanned branches (the number of edges of $\mathbf{T}^b_n$), we have $\mathcal{A}^\star_n = (\mathcal{A}^{i,\star}_n)_{i=1,\dots,K}$ where the r.v. $\mathcal{A}^{i,\star}_n$s are independent with respective distribution $\mathbb{Q}_{\mathcal{H}^i_n}$. The random variable $\Theta^\star_n = (\Theta^\star_n(i))_{i=1,\dots,\kappa}$ has $\kappa$ coordinates that are independent of $(\mathcal{H}_n, \mathcal{A}^\star_n, \mathbf{T}^b_n)$ and distributed as follows:

$$\mathbb{P}(\Theta^\star_n(1) = (j^0, j^1)) = \mu_{j^0} \qquad \text{for any } (j^0, j^1) \in \mathcal{D}_2,$$

$$\mathbb{P}(\Theta^\star_n(i) = (j^0, j^1, j^2)) = 2\mu_{j^0}/\sigma^2_\mu \qquad \text{for any } (j^0, j^1, j^2) \in \mathcal{D}_3, i \geq 2.$$

Since the mean and the variance under $\mu$ are respectively 1 and $\sigma^2_\mu$, these formulas define indeed two distributions. In the following $O^{[l]}$ stands for the vector $(O_1,\dots,O_l)$. Let

$$\Gamma_{n,M} = \{(\mathsf{a}^{[2\kappa-1]},\mathsf{x}^{[2\kappa-1]},\theta^{[\kappa]},T^b),\mathsf{x}_i\in\sqrt{n}[M^{-1},M],\mathsf{a}_i\in J_{\mathsf{x}_i},T^b\in\mathcal{T}^B_{2\kappa-1}\}$$
$$\cap\,\mathrm{supp}(\mathcal{A}^\star_n,\mathcal{H}_n,\Theta^\star_n,\mathbf{T}^b_n).$$

The following Proposition generalizes to finite-dimensional distributions the Proposition 5.

PROPOSITION 15. (i) *For any $\varepsilon>0$, there exists $M$ such that, for $n$ large enough,*

$$\mathbb{P}_n((\mathcal{A}^\star_n,\mathcal{H}_n,\Theta^\star_n,\mathbf{T}^b_n)\in\Gamma_{n,M})\ge 1-\varepsilon.$$

(ii)

$$\sup_{(\mathsf{a},\mathsf{x},\theta,\tau_b)\in\Gamma_{n,M}}\left|\frac{\mathbb{P}_n((\mathcal{A}_n,\mathcal{H}_n,\Theta_n,\mathbf{T}^b_n)=(\mathsf{a},\mathsf{x},\theta,\tau^b))}{\mathbb{P}_n((\mathcal{A}^\star_n,\mathcal{H}_n,\Theta^\star_n,\mathbf{T}^b_n)=(\mathsf{a},\mathsf{x},\theta,\tau^b))}-1\right|\underset{n}{\to}0. \tag{27}$$

In general, $\mathbb{P}_{(\mathcal{A}_n,\mathcal{H}_n,\Theta_n,\mathbf{T}^b_n)}\not\prec\mathbb{P}_{(\mathcal{A}^\star_n,\mathcal{H}_n,\Theta^\star_n,\mathbf{T}^b_n)}$, since $\mathrm{supp}(\Theta^\star_n)$ is strictly included in $\mathrm{supp}(\Theta_n)$ when $\mu[3,+\infty)>0$ (the variable $\Theta^\star$ mimics the coding of binary branchings on $\mathbf{T}^b_n$). In that case, moreover, $\mathbb{P}_n(\mathbf{T}^b_n\notin T^B_{2\kappa-1})>0$ for $n$ large enough, and no control of $(\mathcal{A}_n,\mathcal{H}_n,\Theta_n,\mathbf{T}^b_n)$ is provided on $\complement T^B_{2\kappa-1}$. But an analogous of Lemmas 10 and 11 can be written by weakening slightly the condition $\mathbb{P}_{X_n}\prec\mathbb{P}_{Y_n}$ of Definition 1:

DEFINITION 3. Let $(Y_1,Y_2,\dots)$ and $(X_1,X_2,\dots)$ be two sequences of r.v. taking their values in a Polish space $S$. We say that $\mathbb{P}_{X_n}/\mathbb{P}_{Y_n}\overset{\star}{\to}1$ or $X_n/\!/_\star Y_n\to 1$ if for any $\varepsilon>0$ there exists a measurable set $A^\varepsilon_n$ and a measurable function $f^\varepsilon_n: A^\varepsilon_n\mapsto\mathbb{R}$ satisfying $\mathbb{P}_{X_n}=f^\varepsilon_n\mathbb{P}_{Y_n}$ on $A^\varepsilon_n$, such that $\sup_{x\in A^\varepsilon_n}|f^\varepsilon_n(x)-1|\underset{n}{\to}0$ and such that $\mathbb{P}_{Y_n}(A^\varepsilon_n)\ge 1-\varepsilon$ for $n$ large enough.

LEMMA 16. *Assume that $X_n/\!/_\star Y_n\to 1$, then:*

- *If $Y_n\overset{(d)}{\underset{n}{\to}}Y$, then $X_n\overset{(d)}{\underset{n}{\to}}Y$.*
- *Let $(g_n)$ be a sequence of continuous functions from $S$ into a Polish space $S'$. If $X_n/\!/_\star Y_n\to 1$, then $g_n(X_n)/\!/_\star g_n(Y_n)\to 1$.*

The proof is the same as those of Lemmas 10 and 11 and is left to the reader.

PROOF OF PROPOSITION 15. (i) Let $\varepsilon>0$ be fixed. By Lemma 14, there exists $M$ such that $\mathbb{P}(\Delta_{n,M})>1-\varepsilon/2$. Now conditionally on $\mathbf{T}\in\Delta_{n,M}$, the $\mathcal{H}^i_n$'s belongs to $[\sqrt{n}/M,\sqrt{n}M]$, and then the multinomial random variables $\mathcal{A}_n^{i,\star}$ by Lemma 9 satisfy together $\mathcal{A}_n^{i,\star} \in J_h(\mathcal{H}_n^i)$ with probability arbitrary close to 1 (for any $M$, when $n$ is large enough).

We examine now (ii). Let $(\mathsf{a},\mathsf{x},\theta,\tau^b) \in \Gamma_{n,M}$ for $\theta = ((j_1^0,j_1^1),(j_2^0,j_2^1,j_2^2), \dots,(j_\kappa^0,j_\kappa^1,j_\kappa^2))$.

Using (24) and Comment 1,

$$\begin{aligned}&\mathbb{P}((\mathcal{A}_n,\mathcal{H}_n,\Theta_n,\mathbf{T}_n^b) = (\mathsf{a},\mathsf{x},\theta,\tau^b))\\ &\qquad = \frac{(\prod_{i=1}^{2\kappa-1}\mathbb{Q}_{\mathsf{x}_i}(\mathsf{a}_i))(\prod_{i=1}^{\kappa}\mu_{j_i^0})\mathbb{P}(|\mathbf{f}^{(l)}_{\#\mathcal{S}(\tau^b,\varnothing,\mathsf{a},\theta)}| = F_l(\mathsf{x},\tau^b), 0\le l\le\kappa)}{\mathbb{P}(|\mathbf{T}|=n)},\end{aligned} \tag{28}$$

where the $\mathbf{f}^{(l)}$'s are independent forests.

On the other hand,

$$\begin{aligned}&\mathbb{P}((\mathcal{A}_n^\star,\mathcal{H}_n,\Theta_n^\star,\mathbf{T}_n^b) = (\mathsf{a},\mathsf{x},\theta,\tau^b))\\ &\quad = \mathbb{P}((\mathcal{A}_n^\star,\Theta_n^\star) = (\mathsf{a},\theta)|(\mathcal{H}_n,\mathbf{T}_n^b) = (\mathsf{x},\tau^b))\mathbb{P}((\mathcal{H}_n,\mathbf{T}_n^b) = (\mathsf{x},\tau^b))\\ &\quad = \left(\prod_{i=1}^{2\kappa-1}\mathbb{Q}_{\mathsf{x}_i}(\mathsf{a}_i)\right)\left(\prod_{i=1}^{\kappa}\mu_{j_i^0}\right)\left(\frac{2}{\sigma_\mu}\right)^{\kappa-1}\mathbb{P}((\mathcal{H}_n,\mathbf{T}_n^b) = (\mathsf{x},\tau^b));\end{aligned}$$

summing formula (28) on all possible values of the $\mathsf{a}_i$'s and the $\theta$'s leads to

$$\begin{aligned}&\mathbb{P}((\mathcal{H}_n,\mathbf{T}_n^b) = (\mathsf{x},\tau^b))\\ &\qquad = \frac{(\sigma_\mu^2/2)^{\kappa-1}\mathbb{P}(|\mathbf{f}^{(l)}_{\#\mathcal{S}(\tau^b,\varnothing,\mathsf{m},\tilde\theta)}| = F_l(\mathsf{x},\tau^b), 0\le l\le\kappa)}{\mathbb{P}(|\mathbf{T}|=n)},\end{aligned} \tag{29}$$

where $\mathsf{m} = (\mathsf{m}_i)_{i=1,\dots,\kappa}$ is a vector of $\kappa$ multinomial independent r.v. (the parameters of $\mathsf{m}_i$ are $\mathsf{x}_i$ and $\mathbf{p}$), and where $\tilde\theta = (\tilde\theta(i))_{i\in[\![1,\kappa]\!]} =^{(d)} \Theta_n^\star$ and is independent of $\mathsf{m}$. Hence, for $(\mathsf{a},\mathsf{x},\theta,\tau^b) \in \Gamma_{n,M}$,

$$\begin{aligned}&\frac{\mathbb{P}((\mathcal{A}_n,\mathcal{H}_n,\Theta_n,\mathbf{T}_n^b) = (\mathsf{a},\mathsf{x},\theta,\tau^b))}{\mathbb{P}((\mathcal{A}_n^\star,\mathcal{H}_n,\Theta_n^\star,\mathbf{T}_n^b) = (\mathsf{a},\mathsf{x},\theta,\tau^b))}\\ &\qquad = \frac{\mathbb{P}(|\mathbf{f}^{(l)}_{\#\mathcal{S}(\tau^b,\varnothing,\mathsf{a},\theta)}| = F_l(\mathsf{x},\tau^b), 0\le l\le\kappa)}{\mathbb{P}(|\mathbf{f}^{(l)}_{\#\mathcal{S}(\tau^b,\varnothing,\mathsf{m},\tilde\theta)}| = F_l(\mathsf{x},\tau^b), 0\le l\le\kappa)}.\end{aligned} \tag{30}$$

It is easy to check that for any $(\mathsf{a},\mathsf{x},\tau^b,\theta)$ in $\Gamma_{n,M}$, any $l$, for $n$ large enough,

$$|F_l(\mathsf{x},\tau^b) - n(s_{l+1}-s_l)| \le n^{2/3}, \qquad \left|\#\mathcal{S}_l(\tau^b,\varnothing,\mathsf{a},\theta) - \frac{\sigma_\mu^2}{2}d(u_l,u_{l+1})\right| \le n^{5/12},$$

since $(1/2)^{2/3} < 5/12$. This allows to approximate, on one hand, $F_l(\mathsf{x},\tau^b)$ by $n(s_{l+1}-s_l)$ and, on the other hand, $\#\mathcal{S}_l(\tau^b,\mathsf{a},\theta)$ by $\frac{\sigma_\mu^2}{2}d(u_l,u_{l+1})$ on $\Gamma_{n,M}$ since $n^{5/12} = o(n^{1/2})$, the order of $d(u_l, u_{l+1})$. By Otter and the central local limit theorem and thanks also to a decomposition of the denominator along $\{\mathsf{m} \in \prod J_{\mathsf{x}_i}\}$ or in its complements (as in the proof of Proposition 12), we get

$$\sup_{(\mathsf{a},\mathsf{x},\theta,\tau_b)\in\Gamma_{n,M}} \left| \frac{\mathbb{P}_n(|\mathbf{f}^{(l)}_{\#\mathcal{S}(\tau^b,\varnothing,\mathsf{a},\theta)}| = F_l(\mathsf{x},\tau^b), 0 \le l \le \kappa)}{\mathbb{P}_n(|\mathbf{f}^{(l)}_{\#\mathcal{S}(\tau^b,\varnothing,\mathsf{m},\tilde{\theta})}| = F_l(\mathsf{x},\tau^b), 0 \le l \le \kappa)} - 1 \right| \underset{n}{\to} 0. \tag{31}$$

□

2.3.3. *Proof of the convergence of the finite-dimensional distribution in Theorem 2.* We now show that Proposition 15 and Lemma 16 imply the convergence of the finite-dimensional distributions in Theorem 2. The proof is similar to that of Corollary 13.

Thanks to the Skorohod representation theorem ([14], Theorem 3.30), there exists a probability space $\Omega$ on which the convergence of $\mathbf{h}_n$ to $\mathbf{h}$ is a.s. On $\Omega$, the vector

$$V_n = (\mathbf{h}_n(s_1), \check{\mathbf{h}}_n(s_1,s_2), \mathbf{h}_n(s_2), \check{\mathbf{h}}_n(s_2,s_3), \dots, \mathbf{h}_n(s_\kappa)),$$

which determines $\mathbf{T}_n^b$, as well as the length of the spanned branches, converges a.s. to

$$V_\infty = (\mathbf{h}(s_1), \check{\mathbf{h}}(s_1,s_2), \mathbf{h}(s_2), \check{\mathbf{h}}(s_2,s_3), \dots, \mathbf{h}(s_\kappa)),$$

which determines $\tau_{\mathbf{s}}$ the ordered discrete subtree of the continuum random tree $\tau_\infty$, with contour process $\mathbf{h}$, spanned by the root and the nodes visited at times $s_1, \dots, s_\kappa$. The edge lengths of $\tau_{\mathbf{s}}$ are given by the normalized Brownian excursion (see Aldous [1, 2]). The coordinates of $V_\infty$ are distinct and nonzero a.s., and then $\tau_{\mathbf{s}}$ has only binary branching points, and its shape $\tau_{\mathbf{s}}^b$ belongs to $\mathcal{T}^B_{2\kappa-1}$ (we call here shape the tree $\tau_{\mathbf{s}}$ where the edge lengths are somehow fixed to 1). Let $\mathcal{H}_\infty = (\mathcal{H}_{\infty,i})_{i\in[\![1,2\kappa-1]\!]}$ be the lengths of the (sorted) spanned branches in $\tau_\infty$. By the property of the Brownian excursion, a.s. the coordinate of $\mathcal{H}_\infty$ are almost surely all positive and finite, and then there exists $M$ such that all the $\mathcal{H}_{\infty,i}$ belongs to $[M^{-1}, M]$ (for a $M$ depending on $\tau_\infty$). On $\Omega$,

$$\mathbf{T}_n^b \overset{\text{a.s.}}{\underset{n}{\to}} \tau_{\mathbf{s}}^b, \tag{32}$$

and therefore, for $n$ large enough, $\mathbf{T}_n \in \Gamma_{n,2M}$.

Denote by $(A_n^i)_{i\in[\![1,2\kappa-1]\!]}$ the (sorted) corresponding contents of the spanned branches of $\mathbf{T}_n$, and by $(\mathcal{H}_n^i)_{i\in[\![1,2\kappa-1]\!]}$ their lengths. The normalized contents are then given by

$$\mathbf{g}^{(n)}_{i,k,j} = n^{-1/4}((A_n^i)_{k,j} - \mu_k \mathcal{H}_n^i).$$

Proposition 15 and Lemma 16 entail that

$$((\mathbf{g}_i^{(n)})_{i\in[\![1,2\kappa-1]\!]}, \mathcal{H}_n/\sqrt{n}) \,//\, ((\mathcal{G}^{(i)}(n, \mathcal{H}_n^i))_{i\in[\![1,2\kappa-1]\!]}, \mathcal{H}_n/\sqrt{n}) \to 1,$$

where the r.v. $\mathcal{G}^{(i)}(n,\mathcal{H}_n^i)$'s are independent, and conditionally on $\mathcal{H}_n^i = l$, $(\mathcal{G}^{(i)}(n,\mathcal{H}_n^i)) \overset{(d)}{=} \mathcal{G}(n,l)$ (these variables were introduced in Section 2.1). On $\Omega$, $\mathcal{H}_n/\sqrt{n} \xrightarrow[n]{\text{a.s.}} \mathcal{H}_\infty$. Conditionally on $\mathcal{H}_\infty$, by Proposition 3, $(\mathcal{G}^{(i)}(n, \mathcal{H}_n^i))_{i\in[\![1,2\kappa-1]\!]}$ converges in distribution to $(\mathcal{G}_\infty^{(i)})_{i\in[\![1,2\kappa-1]\!]}$, where the $\mathcal{G}_\infty^{(i)}$ are independent, and $\mathcal{G}_\infty^{(i)}$ is a centered Gaussian vector having as covariance function

$$\operatorname{cov}((\mathcal{G}_\infty^{(i)})_{k,j}, (\mathcal{G}_\infty^{(i)})_{k',j'}) = (-\mu_k\mu_{k'} + \mu_k \mathbb{1}_{(k,j)=(k',j')})\mathcal{H}_{\infty,i}. \tag{33}$$

In order to check that this implies the convergence of the finite-dimensional distributions in Theorem 2, it suffices to reconstitute the contents of the branches $[\![\varnothing, u_i]\!]$'s by summing the contents of the spanned branches they contain, and to use that asymptotically, conditionally on $\mathcal{H}_\infty$, these contents are independent [and that the shape is fixed by (32) for $n$ large enough]. Hence, by (33), one easily gets the fact that each $\mathbf{G}_{k,j}(s_i)$ is Gaussian with the law described in Theorem 2. Knowing $V_\infty$, the limiting $\mathbf{G}(s_i)$'s are obtained as sums of independent Gaussian vectors. To compute the covariance between $\mathbf{G}_{k,j}(s_i)$ and $\mathbf{G}_{k',j'}(s_{i'})$ (for $s_i < s_{i'}$), we use that the nodes in $[\![\varnothing, \check{u}_{i,i'}[\![$ are the common ancestors of $u_i$ and $u_{i'}$. The contents of the branches $]\!]\check{u}_{i,i'}, u_i[\![$ and $]\!]\check{u}_{i,i'}, u_i'[\![$ are asymptotically independent. By (33), one then checks that, knowing $V_\infty$, the covariance $\operatorname{cov}((\mathcal{G}_\infty^{(i)})_{k,j}, (\mathcal{G}_\infty^{(i')})_{k',j'})$ is ruled by the common ancestors, and then equals $(-\mu_k\mu_{k'} + \mu_k\mathbb{1}_{(k,j)=(k',j')}) \times \min\{\mathbf{h}(s), s\in[s_i, s_{i'}]\}$. □

2.4. *Tightness in Theorem 2.* We only prove the tightness of the family $(\mathbf{G}^{(n)})$, since one already knows that $(\mathbf{h}_n)$ is tight [since $\mathbf{h}_n \xrightarrow[n]{(d)} \mathbf{h}$]. In this section we assume $(\mathrm{H}_1)$ and $(\mathrm{H}_2)$.

We collect in the set $\Omega_n^{\alpha,\delta,\gamma,\rho}$ the trees with $n$ edges having some suitable properties:

$$\begin{aligned}\Omega_n^{\alpha,\delta,\gamma,\rho} = \Big\{ & T\in\mathcal{T}_n, \forall t,s\in[0,1], |\mathbf{h}_n(s)-\mathbf{h}_n(t)| \le \delta|t-s|^\alpha, \\ & \max_l ||u(l+1)| - |u(l)|| \le \rho\log n, \\ & |u(n)| < \rho\log n, \forall (k,j)\in I_K, l\in(0,|u|], \\ & \qquad\qquad |A_{u,l,k,j} - \mu_k l| \le \gamma\sqrt{l\log n}\Big\}.\end{aligned}$$

Lemma 17. *For any $\varepsilon > 0$, $\alpha < 1/2$, there exists $\delta > 0$, $\gamma > 0$, $\rho > 0$, s.t.* $\mathbb{P}(\Omega_n^{\alpha,\delta,\gamma,\rho}) \ge 1-\varepsilon$.

According to Lemmas 7 and 8, and Remark 1, only the condition on the Hölderienity of $H$ has to be checked. We refer to Marckert and Miermont [19], Section 5.2, for a proof of this result.

Let $\varepsilon > 0$ be fixed. Set $\alpha = 2/5$ and choose $\delta > 0$, $\gamma > 0$, $\rho > 0$ s.t. $\mathbb{P}(\Omega_n^{\alpha,\delta,\gamma,\rho}) \geq 1 - \varepsilon$ for $n$ large enough. For these choices, write $\Omega_\varepsilon$ instead of $\Omega_n^{\alpha,\delta,\gamma,\rho}$.

We will establish the following proposition.

Proposition 18. *For any $a > 0$, there exist $\beta > 0$, $c > 0$ s.t. for any sufficiently large $n$,*

$$\mathbb{E}(\|\mathbf{G}_s^{(n)} - \mathbf{G}_t^{(n)}\|_1^\beta \mathbb{1}_{\Omega_\varepsilon}) \leq c|t-s|^{1+a} \qquad \text{for any } s,t \in [0,1]. \tag{34}$$

This implies that for any $a > 0$ the $(1+a)/\beta$-Hölder norm of the family $(\mathbf{G}^{(n)})$ is tight, and then that $(\mathbf{G}^{(n)})$ is tight in $C([0,1])^{\#I_K}$ [recall that $\mathbf{G}_0^{(n)}$ is the null vector of $\mathbb{R}^{\#I_K}$].

We first point out that, using (20), we get that for any $a > 0$ there exists $\beta > 0$ such that, for $n$ large enough, $\mathbb{E}(\|\mathbf{G}^{(n)}(s_n) - \mathbf{G}^{(n)}(s)\|_1^\beta \mathbb{1}_{\Omega_n^\rho}) \leq c(s - s_n)^{1+a}$. Hence, we can restrict ourself to prove (34) only for $s$ and $t$ such that $ns$ and $nt$ are integer (this is classical). From now on, we assume that $s,t$ are in $[0,1]_n := [0,1] \cap \mathbb{N}/n$, and $s \neq t$.

We set $u_1 = u(\lfloor ns \rfloor), u_2 = u(\lfloor nt \rfloor)$, $\check{u}_{1,2}$ their deepest common ancestor and $D_n(s,t) = d(u_1,u_2)$. There exists $\delta' > 0$, such that, for $T \in \Omega_\varepsilon$, any $s,t \in [0,1]_n$, $s \neq t$,

$$\begin{aligned} D_n(s,t) &\leq 2 + H_n(nt) + H_n(ns) - 2 \min_{k \in [ns,nt]} H_n(k) \\ &\leq 2\sqrt{n}\delta|t-s|^{2/5} + 2 \leq \delta'\sqrt{n}|t-s|^{2/5}. \end{aligned}$$

Lemma 19. *For any $\alpha' > 0, a > 0$, there exist $\beta > 0, c > 0$ s.t. for any $s,t \in [0,1]$ such that $|s-t| \leq (\log n)^{-3}$, for $n$ large enough,*

$$\mathbb{E}(\|\mathbf{G}_s^{(n)} - \mathbf{G}_t^{(n)}\|_1^\beta \mathbb{1}_{\Omega_\varepsilon}) \leq c|t-s|^{1+a}. \tag{35}$$

Proof. Let $s,t \in [0,1]_n$, $s \neq t$. We use a deterministic bound valid for all trees $T$ in $\Omega_\varepsilon$. Let $(k,j) \in I_K$ fixed. As in the proof of Proposition 4, it suffices to show that

$$n^{-\beta/4}|A_{u_1,k,j} - \mu_k|u_1| - A_{u_2,k,j} + \mu_k|u_2||^\beta \leq c|s-t|^{1+a}. \tag{36}$$

Passing via $\check{u}_{1,2}$, the left-hand side of (36) is smaller than

$$c_1 n^{-\beta/4}(|A_{u_1,h_1,k,j} - \mu_k|h_1||^\beta + |A_{u_2,h_2k,j} - \mu_k|h_2||^\beta + 2^\beta),$$

where $h_1 := d(u_1, \check{u}_{1,2}) - 1$ and $h_2 := d(u_2, \check{u}_{1,2}) - 1$ (the contribution of $\check{u}_{1,2}$ is bounded by the term 2). Using that $|A_{u,l,k,j} - \mu_k l| \le \gamma\sqrt{l \log n}$ for any $l$ and $l \le D_n(s,t) \le \delta' n^{1/2}|t-s|^{2/5}$, we find

$$n^{-\beta/4}|A_{u_1,k,j} - \mu_k |u_1||^\beta + |A_{u_2,k,j} - \mu_k |u_2||^\beta \le c_2(t-s)^{\beta/5}(\log n)^{\beta/2}$$

and since $|t-s| \le (\log n)^{-3}$, $|t-s|^{\beta/6}(\log n)^{\beta/2} \le 1$ and then $c_2(t-s)^{\beta/5} \times (\log n)^{\beta/2}$ is smaller than $|t-s|^{1+a}$ for $\beta$ and $n$ large enough. $\square$

LEMMA 20. *For any $a > 0$, there exist $\beta > 0$, $c > 0$ s.t. for any $t \in [0,1]$, for any $n$ large enough,*

$$\mathbb{E}(\|\mathbf{G}_t^{(n)}\|_1^\beta \mathbb{1}_{\Omega_\varepsilon}) \le ct^{1+a}. \tag{37}$$

PROOF. First, consider the case $t = 1$. In $\Omega_\varepsilon$, we have $|u(n)| \le \rho \log n$ and then

$$\mathbb{E}(\|\mathbf{G}_1^{(n)}\|_1^\beta \mathbb{1}_{\Omega_\varepsilon}) \le c_1(\rho \log n)^\beta n^{-\beta/4}$$

and this is smaller than $c1^{1+a}$ for any $a > 0$, $c > 0$, $\beta > 0$ for $n$ large enough. By the previous lemma and a simple computation [using that $\max ||u(l+1)| - |u(l)|| \le \rho \log n$], one sees that (37) is true if $t \notin V_n$, where

$$V_n := [(\log n)^{-3}, 1 - (\log n)^{-3}].$$

Assume now that $t \in V_n$. In $\Omega_\varepsilon$, the Hölder property of $\mathbf{h}_n$ and the inequality $|u(n)| \le \rho \log n$ imply that, for $t \in V_n$,

$$|u(\lfloor nt \rfloor)| \le \overline{L_n}(t) := c_2 n^{1/2}[t \wedge (1-t)]^\alpha. \tag{38}$$

For any real number $a$, we denote by $a \cdot \mu$ the vector $(a\mu_k)_{(k,j)\in I_K}$. Using (9) and (13), there exists $c_3 > 0$ such that, for $t \in V_n$ and $n$ large enough,

$$\begin{aligned}
&\mathbb{E}(\|\mathbf{G}_t^{(n)}\|_1^\beta \mathbb{1}_{\Omega_\varepsilon}) \\
&\quad \le c_3 \sum_{h \le \overline{L_n}(t)} \sum_{\mathsf{a} \in \mathbb{N}^I_{[h]}} \mathbb{Q}_h(\mathsf{a}) \frac{\|\mathsf{a} - h \cdot \mu\|_1^\beta}{n^{\beta/4-3/2}} \\
&\qquad \times \mathbb{P}(|\mathbf{f}_{N_1(\mathsf{a})}| = \lfloor nt \rfloor - h, |\mathbf{f}'_{1+N_2(\mathsf{a})}| = n + 1 - \lfloor nt \rfloor)
\end{aligned}$$

and by Otter [23]

$$\begin{aligned}
&c_4 \sum_{h \le \overline{L_n}(t)} \sum_{\mathsf{a} \in \mathbb{N}^I_{[h]}} \mathbb{Q}_h(\mathsf{a}) \frac{\|\mathsf{a} - h.\mu\|_1^\beta}{n^{\beta/4-3/2}} \\
&\qquad \times \frac{N_1(\mathsf{a})(1 + N_2(\mathsf{a}))\mathbb{P}(W_{\lfloor nt \rfloor - h} = N_1(\mathsf{a}))\mathbb{P}(W_{n - \lfloor nt \rfloor + 1} = 1 + N_2(\mathsf{a}))}{(\lfloor nt \rfloor - h)(n - \lfloor nt \rfloor + 1)},
\end{aligned}$$

where $(W_k)$ is the random walk described in the beginning of Section 2.2.1. In order to bound these two last probabilities, we use a classical concentration property valid for any nondegenerate random walk $(W_k)_k$ (trivial consequence of Petrov [24], Theorem 2.22 p. 76): there exists a constant $c_5$ such that, for any $n \geq 0$,

$$\sup_y \mathbb{P}(W_n = y) \leq c_5/\sqrt{n}. \tag{39}$$

Now, for any $\mathsf{a} \in \mathbb{N}^I_{[h]}$, $N_1(\mathsf{a})$ and $N_2(\mathsf{a})$ are smaller than $Kh$, and for any $h \leq \overline{L_n}(t)$, $t \in V_n$ and $n$ large enough, $nt - h \geq nt/2$. We then get

$$\mathbb{E}(\|\mathbf{G}^{(n)}_t\|^\beta_1 \mathbb{1}_{\Omega_\varepsilon}) \leq c_6 \sum_{h \leq \overline{L_n}(t)} \sum_{\mathsf{a} \in \mathbb{N}^I_{[h]}} \frac{\mathbb{Q}_h(\mathsf{a}) \|\mathsf{a} - h \cdot \mu\|^\beta_1 h^2}{n^{\beta/4 - 3/2}(\lfloor nt \rfloor - h)^{3/2}(n - \lfloor nt \rfloor + 1)^{3/2}}.$$

Using Proposition 4, we obtain that, for any $t \in V_n$,

$$\mathbb{E}(\|\mathbf{G}^{(n)}_t\|^\beta_1 \mathbb{1}_{\Omega_\varepsilon}) \leq c_7 \frac{(\overline{L_n}(t))^{\beta/2+3}}{n^{\beta/4+3/2}(t(1-t))^{3/2}} \leq c_8 \frac{(t \wedge (1-t))^{\beta/2+3}}{(t(1-t))^{3/2}}. \qquad \square$$

REMARK 2. The last formula implies that, for any $a > 0$, there exist $\beta > 0$, $c > 0$ s.t. for any $t \in V_n$, for any $n$ large enough,

$$\mathbb{E}(\|\mathbf{G}^{(n)}_t\|^\beta_1 \mathbb{1}_{\Omega_\varepsilon}) \leq c(t \wedge (1-t))^{1+a}. \tag{40}$$

This allows to prove a part of Proposition 18: since $\mathbb{E}(\|\mathbf{G}^{(n)}_t - \mathbf{G}^{(n)}_s\|^\beta_1 \mathbb{1}_{\Omega_\varepsilon}) \leq c\mathbb{E}(\mathbb{1}_{\Omega_\varepsilon}(\|\mathbf{G}^{(n)}_t\|^\beta_1 + \|\mathbf{G}^{(n)}_s\|^\beta_1))$ when $s, t \in V_n$ and $s \leq t$,

– if $s \leq t - s$ [in this case $t \leq 2(t-s)$], then $\mathbb{E}(\|\mathbf{G}^{(n)}_t - \mathbf{G}^{(n)}_s\|^\beta_1 \mathbb{1}_{\Omega_\varepsilon}) \leq c(t-s)^{1+a}$,

– if $1 - t \leq t - s$ [in this case $1 - s \leq 2(t-s)$], then

$$\mathbb{E}(\|\mathbf{G}^{(n)}_t - \mathbf{G}^{(n)}_s\|^\beta_1 \mathbb{1}_{\Omega_\varepsilon}) \leq c((1-t)^{1+a} + (1-s)^{1+a}) \leq c_2(t-s)^{1+a}.$$

Thanks to this remark, only the case $s, t \in V_n$, $s \leq t$, and

$$[s \wedge (1-s)] \geq t - s \quad \text{and} \quad [t \wedge (1-t)] \geq t - s \tag{41}$$

remains to be checked. So assume that $s$ and $t$ satisfy these constraints.

Consider $\mathcal{A}_n = (\mathcal{A}^1_n, \mathcal{A}^2_n, \mathcal{A}^3_n) = (A_{(\check{u}_{1,2}, u_1)}, A_{(\check{u}_{1,2}, u_2)}, A_{(\varnothing, \check{u}_{1,2})})$ the contents of the "three" spanned branches in $\mathbf{T}_{\mathbf{s}^2}$ (some of these spanned branches may be empty). We have

$$\begin{aligned} &\mathbb{E}(\|\mathbf{G}^{(n)}_s - \mathbf{G}^{(n)}_t\|^\beta_1 \mathbb{1}_{\Omega_\varepsilon}) \\ &\quad \leq \sum_{h_1,h_2,h_3} \sum_{\mathsf{a}_1,\mathsf{a}_2,\mathsf{a}_3} \frac{\mathbb{P}_n(\mathcal{A}^i_n = \mathsf{a}_i, i = 1,2,3)[\|\mathsf{a}_1 - h_1 \cdot \mu\|^\beta_1 + \|\mathsf{a}_2 - h_2 \cdot \mu\|^\beta_1]}{n^{\beta/4}}, \end{aligned} \tag{42}$$

where the first sum is taken on $h_1 + h_3 \le \overline{L_n}(s)$, $h_2 + h_3 \le \overline{L_n}(t)$, $h_1 + h_2 \le \overline{D_n}(s,t) := \delta' n^{1/2}|t-s|^\alpha$ where $\overline{L_n}(x)$ is given in (38). By (24), Comment 1 and the Otter formula, $h_1, h_2, h_3, \mathsf{a}_1, \mathsf{a}_2, \mathsf{a}_3$ fixed,

$$\mathbb{P}_n(\mathcal{A}_n^i = \mathsf{a}_i, i = 1,2,3) \le cn^{3/2} \sup_\theta \prod_{i=1}^3 \mathbb{Q}_{h_i}(\mathsf{a}_i) \frac{S_i(\theta)\mathbb{P}(W_{F_i} = S_i(\theta))}{F_i}, \tag{43}$$

where the supremum is taken on $\theta = (\theta_1, \theta_2, \theta_3) \in [\![0, K]\!]^8$, and where $F_1 = ns + 1 - |u(\lfloor ns \rfloor)| - 1$, $F_2 = n(t-s) + 1 - (|u(\lfloor nt \rfloor)| - |\check{u}_{1,2}|)$, $F_3 = n(1-t)+1$, $S_1 = N_1(\mathsf{a}_3) + N_1(\mathsf{a}_1) + \theta_1$, $S_2 = N_2(\mathsf{a}_1) + N_1(\mathsf{a}_2) + \theta_2$, $S_3 = N_2(\mathsf{a}_2) + N_2(\mathsf{a}_3) + \theta_3$.

We plug this bound in (42), and bound the left-hand side using the following ingredient:

– the probabilities in (43) involving the random walks are bounded using (39).

– for $\mathsf{a} \in \mathbb{N}^I_{[h]}$, $N_i(\mathsf{a}) \le Kh$ and then for a constant $c > 0$,

$$S_1 \le K|u(\lfloor ns \rfloor)| + \theta_1 \le c\overline{L_n}(s),$$
$$S_2 \le c|\overline{D_n}(s,t)|,$$
$$S_3 \le K|u(\lfloor nt \rfloor)| + \theta_3 \le c\overline{L_n}(t).$$

The denominators are bounded using $|t-s| \ge (\log n)^{-3}$, $[t \wedge (1-t)] \ge (\log n)^{-3}$, $[s \wedge (1-s)] \ge (\log n)^{-3}$, and then for $n$ large enough,

$$F_1 \ge ns/2, \qquad F_2 \ge n(t-1)/2, \qquad F_3 \ge n(1-t)/2.$$

Finally, we get that the left-hand side of (42) is smaller than

$$c\Biggl\{\overline{L_n}(s)\overline{L_n}(t)\overline{D_n}(s,t) \times \sum_{h_1,h_2,h_3} \sum_{\mathsf{a}_1,\mathsf{a}_2,\mathsf{a}_3} \prod_{i=1}^3 \mathbb{Q}_{h_i}(\mathsf{a}_i)[\|\mathsf{a}_1 - h_1 \cdot \mu\|_1^\beta + \|\mathsf{a}_2 - h_2.\mu\|_1^\beta]\Biggr\} \times \{n^{\beta/4-3/2}[n^3(s \wedge (1-s))(t \wedge (1-t))(t-s)]^{3/2}\}^{-1}.$$

The double sum is smaller than

$$\sum_{h_1,h_2,h_3} h_1^{\beta/2} + h_2^{\beta/2} \le (\overline{D_n}(s,t))^{\beta/2+2}\overline{L_n}(s),$$

this last factor $\overline{L_n}(s)$ being a bound of $h_3$. Finally,

$$\mathbb{E}(\|\mathbf{G}_s^{(n)} - \mathbf{G}_t^{(n)}\|_1^\beta \mathbb{1}_{\Omega_\varepsilon}) \le cn^{3/2-\beta/4} \frac{(\overline{L_n}(s))^2\overline{L_n}(t)(\overline{D_n}(s,t))^{\beta/2+3}}{[n^3(s \wedge (1-s))(t \wedge (1-t))(t-s)]^{3/2}}.$$

By (41), it suffices to take $\beta$ large enough. □

2.5. *Proof of Theorem* 1. Consider the representation of $\ell(u)$ given in (2). For any $s$ such that $ns$ is an integer,

$$\mathbf{r}_n(s) = \mathbf{r}_n^{(1)}(s) + \mathbf{r}_n^{(2)}(s), \tag{44}$$

where

$$\mathbf{r}_n^{(1)}(s) = n^{-1/4} \sum_{(k,j)\in I_K} \sum_{l=1}^{A_{u(ns),k,j}} (Y_{k,j}^{(l)} - m_{k,j}),$$
$$\mathbf{r}_n^{(2)}(s) = n^{-1/4} \sum_{(k,j)\in I_K} (A_{u(ns),k,j} - \mu_k |u(ns)|) m_{k,j} = \langle \mathbf{G}^{(n)}(s), \overrightarrow{m} \rangle,$$

where $\overrightarrow{m} = (m_{k,j})_{(k,j)\in I_K}$ and $\langle a, b\rangle = \sum_{(k,j)\in I_K} a_{k,j} b_{k,j}$. For $s$ in $[i/n, (i+1)/, n]$, $\mathbf{r}_n^{(1)}(s)$ and $\mathbf{r}_n^{(2)}(s)$ are defined by linear interpolation. Since $\mathbf{h}_n \overset{(d)}{\underset{n}{\to}} \mathbf{h}$ in $C[0,1]$, by the Skorohod representation theorem [14], Theorem 3.30, there exists a probability space $\Omega$ on which this convergence is a.s. On this space by Theorem 2, $\mathbf{G}^{(n)}$ converges in distribution in $C([0,1])^{\# I_K}$ to $\mathbf{G}^{\mathbf{h}}$, where $\mathbf{G}^{\mathbf{h}}$ has the distribution of $\mathbf{G}$ knowing $\mathbf{h}$. Now, since the application

$$\begin{aligned} \Psi_{\overrightarrow{m}} : C([0,1])^{\# I_K} &\longrightarrow (C[0,1]), \\ (s \mapsto g(s)) &\longmapsto (s \mapsto \langle g(s), \overrightarrow{m}\rangle) \end{aligned}$$

is continuous, on $\Omega$ we have

$$\langle \mathbf{G}^{(n)}, \overrightarrow{m}\rangle \overset{(d)}{\underset{n}{\to}} \mathbf{r}^{(2)} := \langle \mathbf{G}^{\mathbf{h}}, \overrightarrow{m}\rangle \tag{45}$$

in $C([0,1])$. On $\Omega$, $\mathbf{r}^{(2)}$ is a centered Gaussian process with covariance function

$$\begin{aligned} &\operatorname{cov}(\mathbf{r}^{(2)}(s), \mathbf{r}^{(2)}(t)) \\ &\qquad = \check{\mathbf{h}}(s,t) \sum_{(k,j)\in I_K} \sum_{(k',j')\in I_K} (-\mu_k \mu_{k'} + \mu_k \mathbb{1}_{(k,j)=(k',j')}) m_{k,j} m_{k',j'}. \end{aligned}$$

On $\Omega$ (or on an enlarged space), $\mathbf{r}_n^{(1)}$ is the standard head of a discrete snake associated with independent centered displacements. As shown in [19], under ($H_1$) and ($H_2$),

$$\mathbf{r}_n^{(1)} \overset{(d)}{\underset{n}{\to}} \mathbf{r}^{(1)} \tag{46}$$

in $C([0,1],\mathbb{R})$, where $\mathbf{r}^{(1)}$ given $\mathbf{h}$ is a centered Gaussian process with covariance function

$$\operatorname{cov}(\mathbf{r}^{(1)}(s), \mathbf{r}^{(1)}(t)) = \check{\mathbf{h}}(s,t) \sum_{(k,j)\in I_K} \mu_k \sigma_{k,j}^2.$$

It remains to prove that, given $\mathbf{h}$, the finite-dimensional distributions of $\mathbf{r}^{(1)}$ and $\mathbf{r}^{(2)}$ are independent. We establish the "asymptotic independence" between the two processes $\mathbf{r}_n^{(1)}$ and $\mathbf{r}_n^{(2)}$ knowing $\mathbf{h}$. The arguments are quite straightforward; we just explicit the uni-dimensional case. Let

$$\mathcal{T}_n^{\nu} = \{T \in \mathcal{T}_n, \forall (k,j) \in I_K, u \in T, |A_{u,k,j} - \mu_k |u|| \leq n^{1/4+\nu}\}.$$

According to Lemma 9 in [19], for any $\nu > 0$, $\varepsilon > 0$, if $n$ is large enough, $\mathbb{P}_n(\mathcal{T}_n^{\nu}) \geq 1 - \varepsilon$. Letting $s \in [0,1]$ (such that $ns$ is an integer), one may compare

$$\mathbf{r}'_n(s) = n^{-1/4} \sum_{(k,j)\in I_K} \sum_{l=1}^{\lfloor \mu_k |u(\lfloor ns \rfloor)| - n^{1/4+\nu} \rfloor} (Y_{k,j}^{(l)} - m_{k,j})$$

with $\mathbf{r}_n^{(1)}(s)$, where the same r.v. $Y_{k,j}^{(l)}$ are involved in both $\mathbf{r}'_n$ and $\mathbf{r}_n^{(1)}$. Knowing $|u(\lfloor ns \rfloor)|$, $\mathbf{r}'_n(s)$ is independent of $\mathbf{r}_n^{(2)}(s)$ since $\mathbf{r}'_n(s)$ is a function of the $Y_{k,j}$'s when $\mathbf{r}_n^{(2)}$ is a function of the $A_{u(ns),k,j}$'s. We will prove that $|\mathbf{r}_n^{(1)}(s) - \mathbf{r}'_n(s)| \xrightarrow[n]{\text{proba.}} 0$ which is sufficient to deduce that $\mathbf{r}^{(1)}$ and $\mathbf{r}^{(2)}$ are independent given $\mathbf{h}$ (in the uni-dimensional case): indeed, the distance in $\mathbb{R}^3$ between $(\mathbf{r}'_n(s), \mathbf{r}_n^{(2)}(s), \mathbf{h}_n(s))$ and $(\mathbf{r}_n^{(1)}(s), \mathbf{r}_n^{(2)}(s), \mathbf{h}_n(s))$ goes to 0 in probability; hence, $(\mathbf{r}'_n(s), \mathbf{r}_n^{(2)}(s), \mathbf{h}_n(s)) \xrightarrow[n]{(d)} (\mathbf{r}^{(1)}(s), \mathbf{r}^{(2)}(s), \mathbf{h}(s))$ and then $(\mathbf{r}^{(1)}(s), \mathbf{r}^{(2)}(s))$ are independent given $\mathbf{h}(s)$ since $(\mathbf{r}'_n(s), \mathbf{r}_n^{(2)}(s))$ are independent given $(\mathbf{h}_n(s))$.

We have

$$\mathbb{P}_n(|\mathbf{r}'_n(s) - \mathbf{r}_n^{(1)}| \geq x) \leq \mathbb{P}(|\mathbf{r}'_n(s) - \mathbf{r}_n^{(1)}| \geq x, \mathcal{T}_n^{\nu}) + \mathbb{P}_n(\mathcal{T}_n \setminus \mathcal{T}_n^{\nu}).$$

The last term goes to 0 for any $\nu > 0$. The Rosenthal inequality ([24], Theorem 2.11) asserts that if $(X_k)_k$ is a sequence of centered r.v. and $q \geq 2$, then

$$\mathbb{E}\left(\left|\sum_{i=1}^n X_i\right|^q\right) \leq c(q)\left(\sum_{i=1}^n \mathbb{E}(|X_i|^q) + \left(\sum_{i=1}^n \operatorname{var}(X_i)\right)^{q/2}\right), \tag{47}$$

where $c(q)$ is a positive constant depending only on $q$. For $p$ satisfying $(\mathrm{H}_2)$, we have

$$\mathbb{P}(|\mathbf{r}'_n(s) - \mathbf{r}_n^{(1)}| \geq x, \mathcal{T}_n^{\nu}) \leq \mathbb{E}(x^{-p} |\mathbf{r}'_n(s) - \mathbf{r}_n^{(1)}|^p \mathbb{1}_{\mathcal{T}_n^{\nu}}).$$

Conditioning at first by the $A_{(u(ns))}$, and using (47), we get $\mathbb{P}(|\mathbf{r}'_n(s) - \mathbf{r}_n^{(1)}| \geq x, \mathcal{T}_n^{\nu}) \leq$

$$\frac{x^{-p} c(p)}{n^{p/4}} \left(\sum_{(k,j)\in I_K} 2n^{1/4+\nu} \mathbb{E}(|Y_{k,j} - m_{k,j}|^p) + \left(\sum_{(k,j)\in I_K} 2n^{1/4+\nu} \sigma_{k,j}^2\right)^{p/2}\right)$$

and then for $\nu < 1/4$, for any $x > 0$, the bound goes to 0.

Hence, $\mathbf{r}$ is a centered Gaussian process with covariance function sum of the ones of $\mathbf{r}^{(1)}$ and $\mathbf{r}^{(2)}$. Using that $\sum_{(k,j)} \mu_k m_{k,j} = \mathbf{m} = 0$, we get $\operatorname{cov}(\mathbf{r}(s), \mathbf{r}(t)) = \check{\mathbf{h}}(s,t) \sum_{(k,j)} \mu_k \mathbb{E}(Y_{k,j}^2)$.

**Acknowledgments.** I warmly thank Svante Janson who pointed out some errors and imprecisions in a preliminary version of this paper and who made numerous remarks and suggestions. I also thank the referee who helped me to improve the paper.

CNRS, LaBRI
Université Bordeaux 1
351 Cours de la Libération
33405 Talence Cedex
France
E-mail: marckert@labri.fr